\documentclass[preprint,1p,times,lefttitle]{elsarticle}
\journal{arvix}

\biboptions{numbers,sort&compress}
\usepackage{amsmath}
\usepackage{mathtools}
\DeclarePairedDelimiter\ceil{\lceil}{\rceil}
\DeclarePairedDelimiter\floor{\lfloor}{\rfloor}
\usepackage{algorithmicx,algorithm}
\usepackage[noend]{algpseudocode}
\usepackage{caption} 
\usepackage{multirow}
\usepackage{color}
\usepackage{graphicx}
\usepackage{float}
\usepackage{subcaption}
\usepackage{tikz}
\usetikzlibrary{shapes.geometric, arrows}
\tikzstyle{startstop} = [rectangle, rounded corners, minimum height=0.8cm,text centered, draw=black, fill=red!20]
\tikzstyle{arrow} = [thick,->,>=stealth]

\begin{document}

\begin{frontmatter}

\title{Three New Arnoldi-Type Methods for the Quadratic Eigenvalue Problem in Rotor Dynamics\tnoteref{t1}}
\tnotetext[t1]{This research did not receive any specific grant from funding agencies in the public, commercial, or not-for-profit sectors.}

\author[1,2]{LI Dong}
\ead{leslie_ld@126.com}
\author[1,2,3]{CHEN Li-fang\corref{cor1}}
\ead{chenlf@mail.buct.edu.cn}

\cortext[cor1]{corresponding author}

\affiliation[1]{organization={Beijing University of Chemical Technology, College of Mechanical and Electrical Engineering}, postcode={100029}, city={Beijing}, country={China}} 
\affiliation[2]{organization={Key Laboratory of Engine Health Monitoring-Control and Networking (Ministry of Education)}, postcode={100029}, city={Beijing}, country={China}} 
\affiliation[3]{organization={Beijing Key Laboratory of Health Monitoring and self-Recovery for High-End Mechanical Equipment}, postcode={100029}, city={Beijing}, country={China}}

\begin{abstract}
\indent
Three new Arnoldi-type methods are presented to accelerate the modal analysis and critical speed analysis of the damped rotor dynamics finite element (FE) model. 
They are the linearized quadratic eigenvalue problem (QEP) Arnoldi method, the QEP Arnoldi method, and the truncated generalized standard eigenvalue problem (SEP) Arnoldi method. And, they correspond to three reduction subspaces, including the linearized QEP Krylov subspace, the QEP Krylov subspace, and the truncated generalized SEP Krylov subspace, where the first subspace is also used in the existing Arnoldi-type methods.
The numerical examples constructed by a turbofan engine low-pressure (LP) rotor demonstrate that our proposed three Arnoldi-type methods are more accurate than the existing Arnoldi-type methods.
\end{abstract}


\end{frontmatter}

\section{Introduction}\label{section:sec1}
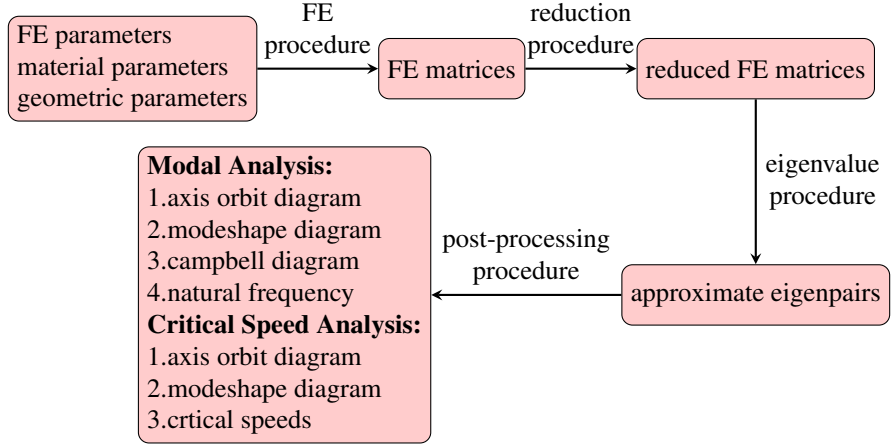
\begin{figure}[ht]
\centering
\begin{tikzpicture}[node distance=2cm]
\node (start) [startstop, align=left] {FE parameters\\ material parameters\\  geometric parameters};
\node (s1) [startstop, right of=start, xshift=2.2cm, align=center] {FE matrices};
\node (s2) [startstop, right of=s1, xshift=2cm, align=center] {reduced FE matrices};
\node (s3) [startstop, below of=s2, yshift=-1cm, align=center] {approximate eigenpairs};
\node (s4) [startstop, left of=s3, xshift=-4.2cm, align=left] {\textbf{Modal Analysis:} \\ 1.axis orbit diagram \\ 2.modeshape diagram \\ 3.campbell diagram \\ 4.natural frequency \\ \textbf{Critical Speed Analysis:} \\ 1.axis orbit diagram \\ 2.modeshape diagram \\ 3.crtical speeds};
\draw [arrow] (start) --node[anchor=south, align=center]{FE \\ procedure} (s1);
\draw [arrow] (s1) --node[anchor=south, align=center]{reduction \\ procedure} (s2);
\draw [arrow] (s2) -- node[anchor=west, align=center]{eigenvalue \\ procedure}(s3);
\draw [arrow] (s3) -- node[anchor=south, align=center]{post-processing \\ procedure}(s4);
\end{tikzpicture}
\caption{The entire process of the modal analysis and critical speed analysis}
\label{fig:FE-Process}
\end{figure}

The modal analysis and critical speed analysis of rotating machinery are performed by a set of procedures, including the FE procedure, eigenvalue procedure, and post-processing procedure. 
First, input FE parameters, material parameters, and geometric parameters, and then generate the FE model with $n$ degrees of freedom or $n$-dimensional FE matrices $M,C,K$ by the FE procedure.
Second, for an undamped rotor, we need to use the eigenvalue procedure to solve the eigenpairs of a $n$-dimensional matrix; for a damped rotor, we need to use the eigenvalue procedure to solve the eigenpairs of a $2n$-dimensional matrix (will be explained in Section~\ref{section:sec2}).
Finally, we can use the post-processing procedure to process eigenpairs, and then output the natural frequency and critical speed and plot the Campbell diagram, axis orbit diagram, and modeshape diagram.
However, for the large-scale and complicated rotating machinery, the FE model has a large number of degrees of freedom, which makes solving the eigenvalue problem by eigenvalue procedure \cite{kublanovskaya1962some, moler1973algorithm, bai2000templates} time-consuming. In rotor dynamics, we are often interested in the first $k$ eigenpairs only, in which eigenpairs are arranged according to the magnitude of eigenvalues from the smallest to the largest. Since $k<<n$, we can reduce the eigenvalue problem first and then solve the reduction problem (shown in Fig.~\ref{fig:FE-Process}).
The reduction methods for the modal analysis and critical speed analysis can be mainly divided into two types: physical subspace methods and Krylov subspace methods \cite{koutsovasilis2008comparison, wagner2010model, sulitka2014application}.
Three new Arnoldi-type methods we proposed belong to the second type.

The first type is proposed from an engineering perspective, and the FE matrices are reduced by being projected into the physical subspace. The methods include static condensation or Guyan reduction \cite{guyan1965reduction}, dynamic condensation \cite{paz1984dynamic}, improved reduced system \cite{gordis1992analysis}, and their variants \cite{friswell1995model, dewen1996improved, QuZuqing1998, friswell2010dynamics}. The advantage of this type is that the elements of the reduced matrices can correspond to the specific node of the FE model, while the disadvantage is that the degrees of freedom need to be partitioned.

The second type is proposed from a mathematical perspective, and the FE matrices are reduced by being projected into the Krylov subspace. The methods can be classified into orthogonal, oblique, and refined projections \cite{lehoucq1996evaluation, bai2000templates, tisseur2001quadratic, bergamaschi2002numerical, guo2004numerical, hernandez2009survey, ozturkmen2022investigation}. When the first two projections use refined eigenvectors to replace Ritz eigenvectors, they can be called refined projections. Advantageously, the degrees of freedom do not need to be partitioned in this type, but the disadvantage is that the element of the reduced matrices is not related to the node of the FE model.

Both the two types have advantages and disadvantages. 
The Krylov subspace methods are superior to the physical subspace methods when considering the dependence on degrees of freedom. Actually, the first type is used in the model updating to compare the FE and experimental modeshape; the second type is used to accelerate the FE analysis.

\begin{figure}[ht]
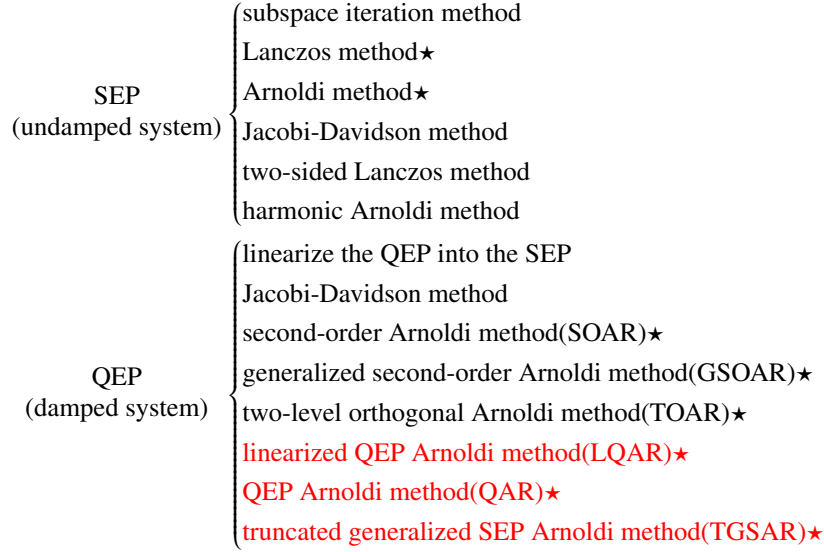

\centering
\begin{align*}
&\parbox{2.8cm}{\centering SEP\\ (undamped system)}
	\begin{cases*}
	\text{subspace iteration method}\\
	\text{Lanczos method$\star$}\\
	\text{Arnoldi method$\star$}\\
	\text{Jacobi-Davidson method}\\
	\text{two-sided Lanczos method}\\
	\text{harmonic Arnoldi method}\\
	\end{cases*}\\
&\parbox{2.8cm}{\centering QEP\\ (damped system)}
	\begin{cases*}
	\text{linearize the QEP into the SEP}\\
	\text{Jacobi-Davidson method}\\
	\text{second-order Arnoldi method(SOAR)$\star$}\\
	\text{generalized second-order Arnoldi method(GSOAR)$\star$}\\
	\text{two-level orthogonal Arnoldi method(TOAR)$\star$}\\
	\text{\textcolor{red}{linearized QEP Arnoldi method(LQAR)$\star$}}\\
	\text{\textcolor{red}{QEP Arnoldi method(QAR)$\star$}}\\
	\text{\textcolor{red}{truncated generalized SEP Arnoldi method(TGSAR)$\star$}}
	\end{cases*}
\end{align*}
\vspace*{-5mm}
\caption{The Krylov subspace methods classified by the type of eigenvalue problem}
\label{fig: Krylov subspace methods}
\end{figure}

The Krylov subspace methods are shown in Fig.~\ref{fig: Krylov subspace methods}, which can be used to reduce the SEP and QEP. Considering our needs of computing multiple continuous right eigenpairs at a time with higher accuracy and efficiency, the orthogonal projection Arnoldi-type methods (shown in Fig.~\ref{fig: Krylov subspace methods} with star marker) are better than others. The Lanczos and Arnoldi methods\cite{lanczos1950iteration, arnoldi1951principle} can be well applied to reduce the SEP. 
However, the Arnoldi, SOAR \cite{ bai2004second, bai2005soar}, GSOAR \cite{otto2004arnoldi}, TOAR \cite{su2008compact, lu2016stability, meerbergen2018mixed}, and their variants \cite{jia2015implicitly, gong2016rgsoar, gong2017rtoar, ravibabu2019stability} cannot achieve enough accuracy for reducing the QEP. Therefore, we propose three Arnoldi-type methods (shown in Fig.~\ref{fig: Krylov subspace methods} with red color) to reduce the QEP. 

The rest of this paper organizes as follows. In Section~\ref{section:sec2}, we introduce the eigenvalue problem in rotor dynamics. In Section~\ref{section:sec3},  we present three Arnoldi-type methods and discuss their relationship with the existing Arnoldi-type methods. To evaluate the capability and performance of the proposed methods, numerical examples are given in Section~\ref{section:sec4}.
Concluding remarks are in Section~\ref{section:sec5}.

The descriptions of the main symbols used in the paper are listed in Table~\ref{table:Notation}, and others that are not included will be explained as used.

\begin{table}[H]
\caption{The description of the main symbols} 
\centering
\begin{tabular}{cl}
\hline
\textbf{Symbols}    & \textbf{Description}                              
\\ 
\hline
I                       & identity matrix                                   
\\
0                       & zero vector or matrix                             
\\
$\mathcal{S}_m$   &  $m$-dimensional SEP Krylov subspace, $\mathcal{S}_m(A;b) = span\{b, Ab, ..., A^{m-1}b\}$
\\
$\mathcal{Q}_m^1$         &  $m$-dimensional linearized QEP Krylov subspace  
\\
$\mathcal{Q}_m^2$         &  $m$-dimensional QEP Krylov subspace  
\\
$\mathcal{G}_{2^m-1}$         & $(2^m-1)$-dimensional generalized SEP Krylov subspace             
\\
$\tilde{\mathcal{G}}_{2m-1}$ & $(2m-1)$-dimensional truncated generalized SEP Krylov subspace  
\\
$V_{\mathcal{S}}$                  & orthonormal basis of $\mathcal{S}_m$                \\
$V_{\mathcal{Q}^1}$                  & orthonormal basis of $\mathcal{Q}_m^1$                \\
$V_{\mathcal{Q}^2}$                  & orthonormal basis of $\mathcal{Q}_m^2$                \\
$V_{\mathcal{G}}$                  & orthonormal basis of $\mathcal{G}_{2^m-1}$                \\
$V_{\tilde{\mathcal{G}}}$        & orthonormal basis of $\tilde{\mathcal{G}}_{2m-1}$        \\
$b$,  $b_0$,  $b_1$     & starting vector                       
\\ 
${\cdot}^H$        &  matrix conjugate transpose  
\\ 
$\bar{\cdot}$ & complex conjugate
\\
\hline
\end{tabular}
\label{table:Notation}
\end{table}

\section{The Eigenvalue Problem in Rotor Dynamics}\label{section:sec2}
In the stationary reference frame, the general rotor dynamics equation is \nocite{meirovitch1980computational, adams1981insights, adams1987insights, ZhongYie1987, friswell2010dynamics}: 
\begin{equation}
M^s \ddot{x} + (C^s + C^{ss})\dot{x} + (K^s + K^{ss})x = f
\end{equation}
where $M^s$ is the conservative mass matrix,  $C^s$ is the non-conservative damping matrix, $K^s$ is the conservative stiffness matrix, $C^{ss}$ is the conservative gyroscopic matrix, $K^{ss}$ is the non-conservative circulatory matrix, $f$ is the external force vector, $.^s$ denotes symmetry, $.^{ss}$ denotes skew-symmetry. They are $n$-dimensional square matrices or vectors.

For a linear rotor-support-foundation system, the rotor system matrices can be expressed explicitly as: $M^s = M$, $C^s = C_0$, $C^{ss} = \Omega C_1$, $K^s = K_0$, $K^{ss} = \Omega K_1$, leading to the following equation:
\begin{equation}
M \ddot{x} + C \dot{x} + K x = f 
\label{eq:2}
\end{equation}
where $C = C_0 + \Omega C_1$, $K = K_0 + \Omega K_1$, and $\Omega$ is the angular velocity (rad/s). $M$ is always positive-definite. For the constrained system, $K$ is positive-definite; for the free system, $K$ is positive semi-definite. 

\subsection{Modal Analysis}\label{section:sec2.1}
Let $x=v e^{st}$, $f=0$, where $v$ is a complex vector, $s$ is a complex number,  Eq.~(\ref{eq:2}) can be rearranged as:
\begin{equation}
(s^2 M + s C + K)v = 0
\label{eq:2-0}
\end{equation}

For an undamped system, $C=0$, Eq.~(\ref{eq:2-0}) is the generalized eigenvalue problem (GEP), we can transform it into the SEP:
\begin{equation}
-M^{-1}K v = s^2 v
\label{eq:2-1}
\end{equation}
We can get $n$ eigenpairs $(s_i^2, v_i)$, where $s_i^2 = w_i^2$, $w_i$ is the undamped natural frequency, $v_i$ is the undamped modeshape.

For a damped system, $C \neq 0$, Eq.~(\ref{eq:2-0}) is the QEP, we also can transform it into the SEP:
\begin{equation}
\begin{bmatrix} -M^{-1}C & -M^{-1}K \\ I & 0\end{bmatrix} \begin{bmatrix} s v \\ v \end{bmatrix} = s\begin{bmatrix} s v \\ v \end{bmatrix}
\label{eq:2-2}
\end{equation}
We can get $2n$ conjugate eigenpairs $(s_i, v_i)$, where $s_{n+i} = \bar{s}_i$, $v_{n+i} = \bar{v}_i$, $s_i = -\xi_i w_i + j w_{di}$, $w_{di} = \sqrt{1 - \xi_i^2} w_i$, $w_{di}$ is the damped natural frequency, $w_i$ is the undamped natural frequency, $\xi_i$ is the damping ratio, $v_i$ is the damped modeshape, $._d$ denotes the system has damping. As $\xi_i$ tends to $0$, $w_{di}$ will tend to $w_{i}$.

\subsection{Critical Speed Analysis}\label{section:sec2.2}
Let $x=v e^{j w_f t}$, $f=0$,  where $v$ is a complex vector, $w_f = n \Omega$ is the excitation frequency, $n$ is the ratio of excitation frequency to angular velocity, Eq.~(\ref{eq:2}) can be rearranged as:
\begin{equation}
(\Omega^2 \hat{M} + \Omega \hat{C} + \hat{K})v = 0
\label{eq:2-3}
\end{equation}
where $\hat{M} = -n^2M + jnC_1$, $\hat{C}=jnC_0+K_1$, $\hat{K}=K_0$. Solving the undamped system $\hat{C}=0$ and damped system $\hat{C} \neq 0$ is same as Section~\ref{section:sec2.1}, the eigenvalue $\Omega$ is critical speed, and the eigenvector $v$ is modeshape.  Although $\Omega$ is a real number, the solution is always a complex number, we need to discuss this situation: if $|Im(\Omega)| > |Re(\Omega)|$, it represents this critical speed does not exist in the real world, in the Campbell diagram, we will find that the natural frequency line has no intersect point with the excitation frequency line; if $|Im(\Omega)| < |Re(\Omega)|$, then $Re(\Omega)$ is the critical speed.

In a broad sense, the critical speed analysis is a special case of the modal analysis, for the solution of the former is a proper subset of the latter. When natural frequency $w$ or $w_d$ equals to excitation frequency $w_f$,  angular velocity $\Omega$ is called critical speed.

\section{Three New Arnoldi-type Methods}\label{section:sec3}
In rotor dynamics, we are interested in the first $k$ eigenpairs, but applying the Arnoldi-type methods to reduce the original problem will converge to the $m$ dominant eigenpairs, where $k \leq m \ll n$. Therefore, we should solve the inverse problem.

We use equations of modal analysis to derive, for the equations of modal analysis and critical speed analysis have the same form. The inverse of Eq.~(\ref{eq:2-1}) is given by:
\begin{equation}
-K^{-1}M v = \frac{1}{s^2} v
\label{eq:3-0-1}
\end{equation}
The inverse of Eq.~(\ref{eq:2-2}) is given by:
\begin{equation}
\begin{bmatrix} 0 & I \\ B & A \end{bmatrix} \begin{bmatrix} s v \\ v \end{bmatrix} = \frac{1}{s} \begin{bmatrix} s v \\ v \end{bmatrix}
\label{eq:3-0-2}
\end{equation}
where $A = -K^{-1}C$, $B = -K^{-1}M$.

For solving Eq.~(\ref{eq:3-0-1}), the reduction subspace is $\mathcal{S}_m(-K^{-1}M;b)$, if $-K^{-1}M$ is Hermitian, then we apply the Lanczos method to generate $V_{S}$; if $-K^{-1}M$ is non-Hermitian, then we apply the Arnoldi method to generate $V_{S}$. By imposing the Galerkin condition, $-K^{-1}M$ can be reduced as $V_{S}^H (-K^{-1}M) V_{S}$, then solve the reduced SEP by the eigenvalue procedure.

For solving Eq.~(\ref{eq:3-0-2}), we give three more accurate reduction methods in the following, which we call the LQAR, QAR, and TGSAR methods, and then discuss their relationship with the existing Arnoldi-type methods.

\subsection{Three QEP Reduction Subspaces}\label{section:sec3.1}
Eq.~(\ref{eq:3-0-2}) can be written as:
\begin{equation}
sBv + Av = \frac{1}{s}v
\label{eq:3-1-1}
\end{equation}
where $s \in C$, $v \in C^n$, and $A, B \in C^{n\times n}$. The orthogonal projection Krylov subspace methods are to find an approximate eigenpair $(\tilde{s}, \tilde{v})$ of Eq.~(\ref{eq:3-1-1}) on the QEP Krylov subspace $\mathcal{Q}_m$, and residual vector satisfies the Galerkin condition:
\begin{equation}
(\tilde{s}B\tilde{v} + A\tilde{v} - \frac{1}{\tilde{s}}\tilde{v}) \perp \mathcal{Q}_m
\end{equation}
or, equivalently, 
\begin{equation}
q^H(\tilde{s}B\tilde{v} + A\tilde{v} - \frac{1}{\tilde{s}}\tilde{v})=0, \forall q \in \mathcal{Q}_m
\end{equation}
where $\tilde{s} \in C$ and $\tilde{v} \in \mathcal{Q}_m$. 
Since $\tilde{v} \in \mathcal{Q}_m$, and then $\tilde{v}$ can be linearly represented by the orthonormal basis $V_{\mathcal{Q}}$ of the subspace $\mathcal{Q}_m$:
\begin{equation}
\tilde{v} = V_{\mathcal{Q}}w_{\mathcal{Q}}
\end{equation}
where, $V_{\mathcal{Q}} = [q_1,q_2, ...,q_m]$. The Galerkin condition can be rewritten as:
\begin{equation}
q^H_i(\tilde{s}BV_{\mathcal{Q}}w_{\mathcal{Q}}+ AV_{\mathcal{Q}}w_{\mathcal{Q}} - \frac{1}{\tilde{s}}V_{\mathcal{Q}}w_{\mathcal{Q}})=0, i=1,2,...,m
\end{equation}
Represent the $m$ equations in matrix form:
\begin{equation}
V_{\mathcal{Q}}(\tilde{s}BV_{\mathcal{Q}}w_{\mathcal{Q}}+ AV_{\mathcal{Q}}w_{\mathcal{Q}} - \frac{1}{\tilde{s}}V_{\mathcal{Q}}w_{\mathcal{Q}})=0
\end{equation}
Finally, rewrite the above equation as:
\begin{equation}
\tilde{s}B_{\mathcal{Q}} w_{\mathcal{Q}} + A_{\mathcal{Q}} w_{\mathcal{Q}} = \frac{1}{\tilde{s}} w_{\mathcal{Q}}
\label{eq:3-1-2}
\end{equation}
where $B_{\mathcal{Q}} = V_{\mathcal{Q}}^H B V_{\mathcal{Q}}$, $A_{\mathcal{Q}} = V_{\mathcal{Q}}^H A V_{\mathcal{Q}}$, and $w_{\mathcal{Q}} = V_{\mathcal{Q}}^H \tilde{v}$. 
Therefore, the essence of reducing the QEP is the reduction of the matrices $B$ and $A$. By finding an orthonormal basis $V_{\mathcal{Q}}$ of $\mathcal{Q}_m$, the matrices $B$ and $A$ can be projected into $\mathcal{Q}_m$ to achieve the reduction. 

Eq.~(\ref{eq:3-1-2}) needs to be linearized before solving:
\begin{equation}
\begin{bmatrix} 0 & I \\ B_{\mathcal{Q}} & A_{\mathcal{Q}} \end{bmatrix} \begin{bmatrix} \tilde{s} w_{\mathcal{Q}} \\ w_{\mathcal{Q}} \end{bmatrix} = \frac{1}{\tilde{s}} \begin{bmatrix} \tilde{s} w_{\mathcal{Q}} \\ w_{\mathcal{Q}} \end{bmatrix}
\label{eq:3-1-3}
\end{equation}
where $\begin{bmatrix} 0 & I \\ B_{\mathcal{Q}} & A_{\mathcal{Q}} \end{bmatrix} = \begin{bmatrix} V_{\mathcal{Q}} & 0 \\ 0 & V_{\mathcal{Q}} \end{bmatrix}^H \begin{bmatrix} 0 & I \\ B & A \end{bmatrix} \begin{bmatrix} V_{\mathcal{Q}} & 0 \\ 0 & V_{\mathcal{Q}} \end{bmatrix}$, 
$\begin{bmatrix} \tilde{s} w_{\mathcal{Q}} \\ w_{\mathcal{Q}} \end{bmatrix} = \begin{bmatrix} V_{\mathcal{Q}} & 0 \\ 0 & V_{\mathcal{Q}} \end{bmatrix}^H \begin{bmatrix} \tilde{s} \tilde{v} \\ \tilde{v} \end{bmatrix}$.
Thus, the essence of reducing the QEP can also be seen as the reduction of the matrix $\begin{bmatrix} 0 & I \\ B & A \end{bmatrix}$. Through an orthonormal basis $V_{\mathcal{Q}}$, the matrix $\begin{bmatrix} 0 & I \\ B & A \end{bmatrix}$ can be projected into $span\{ \begin{bmatrix} V_{\mathcal{Q}} & 0 \\ 0 & V_{\mathcal{Q}} \end{bmatrix} \}$ to achieve the reduction. 

We expect the exact eigenvector $v$ of the QEP and $V_{\mathcal{Q}}$ have the relationship as:
\begin{equation}
span\{\begin{bmatrix} v_1, ..., v_m \end{bmatrix}\} = span\{V_{\mathcal{Q}}\}
\end{equation}
\begin{equation}
span\{\begin{bmatrix} s_1 v_1 & \cdots & s_{m} v_{m} & \bar{s}_1 \bar{v}_1 & \cdots & \bar{s}_{m} \bar{v}_{m}\\ v_1 & \cdots & v_{m} & \bar{v}_1 & \cdots & \bar{v}_{m}\end{bmatrix}\} = span\{\begin{bmatrix} V_{\mathcal{Q}} & 0 \\ 0 & V_{\mathcal{Q}} \end{bmatrix}\}
\end{equation}
The more accurate the subspace $\mathcal{Q}_m$ is, the closer approximation eigenpair $(\tilde{s},\tilde{v})$ will be to the exact eigenpair $(s,v)$. Similar like the subspace $\mathcal{S}_m$ can be spanned by the sequences generated in the power iteration process of the SEP, we present two type subspaces $\mathcal{Q}_m^i$ which can be spanned by the sequences generated in the power iteration process of the linearized QEP and the QEP, respectively.

For the linearized QEP Eq.~(\ref{eq:3-0-2}), the power iteration process of the matrix
$\begin{bmatrix} 0 & I \\ B & A \end{bmatrix}$ is as follows:
\begin{align*}
& \begin{bmatrix} r_{-1} \\ r_0 \end{bmatrix} = \begin{bmatrix} b_0 \\ b_1 \end{bmatrix}\\
& \begin{bmatrix} r_0 \\ r_1 \end{bmatrix} = \begin{bmatrix} 0 & I \\ B & A \end{bmatrix} \begin{bmatrix} b_0 \\ b_1 \end{bmatrix} =  \begin{bmatrix} b_1 \\ Bb_0+Ab_1 \end{bmatrix}\\
& \begin{bmatrix} r_1 \\ r_2 \end{bmatrix} = \begin{bmatrix} 0 & I \\ B & A \end{bmatrix} \begin{bmatrix} b_1 \\ Bb_0+Ab_1 \end{bmatrix} = \begin{bmatrix} Bb_0+Ab_1 \\ Bb_1+ABb_0+A^2b_1 \end{bmatrix}\\
& \begin{bmatrix} r_2 \\ r_3 \end{bmatrix} = \begin{bmatrix} 0 & I \\ B & A \end{bmatrix} \begin{bmatrix} Bb_0+Ab_1 \\ Bb_1+ABb_0+A^2b_1 \end{bmatrix}=  \begin{bmatrix} Bb_1+ABb_0+A^2b_1 \\ B^2b_0+BAb_1+ABb_1+A^2Bb_0+A^3b_1 \end{bmatrix}\\
& ...
\end{align*}
The linearized QEP Krylov subspace is defined by the above sequences as:
$$
\mathcal{Q}_m^1 = span\{b_1, Bb_0+Ab_1, Bb_1+ABb_0+A^2b_1, B^2b_0+BAb_1+ABb_1+A^2Bb_0+A^3b_1, ...\}
$$
or, 
$$
\mathcal{Q}_m^1 = span\{r_0,r_1, r_2, ..., r_{m-1}\}
$$
where $r_j = Br_{j-2} + Ar_{j-1}, j \geq 1$ and $r_{-1} = b_0, r_0 = b_1$.
We can use its orthonormal basis $V_{\mathcal{Q}^1}$ to achieve reduction.

For the QEP Eq.~(\ref{eq:3-1-1}), it can be rewritten as:
\begin{equation}
(sB + A)v = \frac{1}{s}v
\label{eq:3-1-8}
\end{equation}
Similar to the linearized QEP Krylov subspace, the QEP Krylov subspace is defined as:
$$
\mathcal{Q}_m^2 = span\{b, (sB + A)b, (sB + A)^2b, (sB + A)^3b, ...\}
$$
We can use its orthonormal basis $V_{\mathcal{Q}^2}$ to achieve reduction.

We can construct $\mathcal{Q}_m^2$ by solving an eigenvalue, but now we give another subspace which can replace $\mathcal{Q}_m^2$ and donot need eigenvalue. Each sequence of $\mathcal{Q}_m^2$ can be represeneted by the linear combination of $2^i$ sequences:
\begin{align*}
& \textcolor{red}{2^0}: b \\
& \textcolor{red}{2^1}: Bb, Ab \\
& \textcolor{red}{2^2}: B^{2} b, A B b, B A b, A^{2} b\\
& \textcolor{red}{2^3}: B^{3} b, B^{2} A b, A B^{2} b, B A B b, A^{2} B b , A B A b, B A^{2} b, A^{3} b\\
& \textcolor{red}{2^4}: B^{4} b, B^{2} A B b, B^{3} A b, A B^{3} b, B A B^{2} b, A^{2} B^{2} b, B^{2} A^{2} b, A B^{2} A b, A B A B b, B A^{2} B b, B A B A b,\\ &A^{2} B A b, A^{3} B b, A B A^{2} b, B A^{3} b, A^{4} b\\
& \textcolor{red}{2^5}: B^{5} b, B^{2} A B^{2} b, B^{3} A B b, B^{4} A b, A B^{4} b, B A B^{3} b, A^{2} B^{3} b, B^{2} A^{2} B b, B^{2} A B A b, B^{3} A^{2} b, A B^{2} A B b,\\ &A B^{3} A b, A B A B^{2} b, B A^{2} B^{2} b, B A B^{2} A b, B A B A B b, B^{2} A b, A^{2} B A B b, A^{3} B^{2} b, B^{2} A^{3} b, A B^{2} A^{2} b,\\ &A B A^{2} B b , A B A B A b, B A^{2} B A b, B A^{3} B b, B A B A^{2} b , A^{2} B A^{2} b, A^{3} B A b, A^{4} B b, A B A^{3} b, B A^{4} b, A^{5} b\\
& ...
\end{align*}
Therefore, we define $\mathcal{G}_{2^m-1} = \mathcal{S}_m(A;b) \cup \mathcal{S}_m(B;b) \cup \mathcal{S}_{\ceil*{\frac{m}{2}}}(AB;b) \cup \mathcal{S}_{\ceil*{\frac{m}{2}}}(BA;b) \cup \mathcal{S}_{\ceil*{\frac{m}{3}}}(A^2B;b) \cup \mathcal{S}_{\ceil*{\frac{m}{3}}}(B^2A;b) \cup \mathcal{S}_{\ceil*{\frac{m}{3}}}(AB^2;b)\cup \mathcal{S}_{\ceil*{\frac{m}{3}}}(BA^2;b) \cup \mathcal{S}_{\ceil*{\frac{m}{3}}}(ABA;b) \cup \mathcal{S}_{\ceil*{\frac{m}{3}}}(BAB;b) \cup ...$, where $\ceil*{\cdot}$ is the ceil function. We call $\mathcal{G}_{2m-1}$ the generalized SEP Krylov subspace.
$\mathcal{Q}_m^2$ and $\mathcal{G}_{2^m-1}$ have the relationship as:
\begin{equation}
\mathcal{Q}_m^2 \subset \mathcal{G}_{2^m-1}
\end{equation}

However, the dimension of $\mathcal{G}_{2^m-1}$ is large, which leads to generating $V_{\mathcal{G}}$ time-consuming. In practice, for the rotor dynamics FE model reduction problem, we will obtain a satisfactory accuracy when only truncate $\mathcal{G}_{2^m-1}$ as:
$$\tilde{\mathcal{G}}_{2m-1} =  \mathcal{S}_m(A;b) \cup \mathcal{S}_m(B;b)
$$
We call $\tilde{\mathcal{G}}_{2m-1}$ the truncated generalized SEP Krylov subspace, and we can use its orthonormal basis $V_{\tilde{\mathcal{G}}}$ to achieve reduction.

So far, we have presented three subspaces in this section, they are $\mathcal{Q}_m^1$, $\mathcal{Q}_m^2$, and $\tilde{\mathcal{G}}_{2m-1}$. 
In the next, we will give three Arnoldi-type procedures to generate the orthonormal basis for these three subspaces, which are the linearized QEP Arnoldi (LQAR) procedure, the QEP Arnoldi (QAR) procedure, and the truncated generalized SEP Arnoldi (TGSAR) procedure.

\subsection{LQAR Procedure for the Subspace $\mathcal{Q}_m^1$}
\begin{algorithm}[ht]
\begin{algorithmic}[1] 
\State $V = [b_1/||b_1||_2]$
\State $r_0 = b_0/||b_0||_2$
\State $r_1 = b_1/||b_1||_2$
\For {$j = 1$ to $m$}
	\State $\hat{r}_0 = r_1$
	\State $r_1 = Br_0 + Ar_1$
	\State $r_0 = \hat{r}_0$
	\State $\alpha = ||r_1||_2$
	\For {$i = 1$ to $j$}
		\State $r_1 = r_1 - (r_1^H V[:,i])V[:,i]$, where $V[:,i]$ means $i$th column of $V$
	\EndFor
	\If {$||r_1||_2 < \eta*\alpha$}
		\For {$i = 1$ to $j$}
			\State $r_1 = r_1 - (r_1^H V[:,i])V[:,i]$
		\EndFor
	\EndIf
	\If {$||r_1||_2 = 0$}
		\State stop
	\Else
		\State $V = [V; r_1/||r_1||_2]$
	\EndIf
\EndFor
\State $V_{\mathcal{Q}^1} = V[:,1:m]$, where $V[:,1:m]$ means from $1$th column to $m$th column of $V$
\end{algorithmic}
\caption{The LQAR procedure with the partial reorthogonalization technique} 
\label{algorithm:Arnoldi1}
\end{algorithm}

Based on the Arnoldi procedure, we give the LQAR procedure in Algorithm~\ref{algorithm:Arnoldi1}. In addition, due to the finite precision computation, as $m$ increases, $V_{\mathcal{Q}^1}$ will not be an orthonormal basis of $\mathcal{Q}^1_m$. Same as the existing Arnoldi-type methods, we can use the partial reorthogonalization technique \cite{parlett1998symmetric, bai2000templates} to ensure that $V_{\mathcal{Q}^1}$ is still an orthonormal basis of $\mathcal{Q}^1_m$,  hence, the partial reorthogonalization technique is added in Algorithm~\ref{algorithm:Arnoldi1}.

\subsection{QAR Procedure for the Subspace $\mathcal{Q}_m^2$}
\begin{algorithm}[ht]
\begin{algorithmic}[1] 
\State $u_0 = b_0/||b_0||_2$
\State $u_1 = b_1/||b_1||_2$
\For {$k = 1$ to $m-1$}
	\State $\hat{u}_0 = u_1$
	\State $u_1 = Bu_0 + Au_1$
	\State $u_0 = \hat{u}_0$
\EndFor
\State $s = ||u_0||/||u_1||$
\State $V = [b/||b||_2]$
\For {$j = 1$ to $m$}
	\State $w=(sB+A) V[:,j]$
	\State $\alpha = ||w||_2$
	\For {$i = 1$ to $j$}
		\State $w = w - (w^H V[:,i])V[:,i]$
	\EndFor
	\If {$||w||_2 < \eta*\alpha$}
		\For {$i = 1$ to $j$}
			\State $w = w - (w^H V[:,i])V[:,i]$
		\EndFor
	\EndIf
	\If {$||w||_2 = 0$}
		\State stop
	\Else
		\State $V = [V; w/||w||_2]$
	\EndIf
\EndFor
\State $V_{\mathcal{Q}^2} = V[:,1:m]$
\end{algorithmic}
\caption{The QAR procedure with the partial reorthogonalization technique} 
\label{algorithm:Arnoldi2}
\end{algorithm}

As we all kown, by power iteration, we will get the  approximate dominant eigenvector $\begin{bmatrix} \tilde{s}_1 \tilde{v} \\ \tilde{v} \end{bmatrix}$ of the matrix $\begin{bmatrix} 0 & I \\ B & A \end{bmatrix}$.
Thus, we can use $\tilde{s}_1$ to construct $\mathcal{Q}_m^2$, where $\tilde{s}_1 = ||\tilde{s}_1 \tilde{v}|| / ||\tilde{v}||$. And, the orthonormal basis $V_{\mathcal{Q}^2}$ of $\mathcal{Q}_m^2$ can be generated by the Arnoldi procedure. Consequently, we give the QAR procedure with the partial reorthogonalization technique in Algorithm~\ref{algorithm:Arnoldi2}.

\subsection{TGSAR Procedure for the Subspace $\tilde{\mathcal{G}}_{2m-1}$}
We give the first version of the TGSAR procedure in Algorithm~\ref{algorithm:Arnoldi3} to generate the orthonormal basis of $\tilde{\mathcal{G}}_{2m-1}$. The computational time is mainly related to the number of the modified Gram-Schmidt iterations.
The first version includes three orthogonalization processes, which generate orthonormal basis $V_1, V_2, V_{\tilde{\mathcal{G}}}$ respectively. 
Three processes total need $N_1 = \frac{m(m+1)}{2} + \frac{m(m+1)}{2} + \frac{3m(m-1)}{2} = \frac{5m^2 - m}{2}$ iterations. 
We find that three processes can be merged into one process, and then only requires $N_2=\frac{2m(2m-1)}{2}= 2m^2 - m$ iterations.
Therefore, we give the second version of the TGSAR procedure with the partial reorthogonalization technique in Algorithm~\ref{algorithm:Arnoldi4}.

\begin{algorithm}[ht]
\begin{algorithmic}[1] 
\State $V_1 = [b/||b||_2], V_2=[b/||b||_2]$
\For {$j = 1$ to $m$}
\State $w = A V_1[:,j]$
	\For {$i = 1$ to $j$}
	\State $w = w - (w^H V_1[:,i]) V_1[:,i]$
	\EndFor
	\If {$||w|| = 0$}
		\State stop
	\Else
		\State $V_1 = [V_1; w/||w||_2]$
	\EndIf
\EndFor
\For {$j = 1$ to $m$}
\State $w = B V_2[:,j]$
	\For {$i = 1$ to $j$}
	\State $w = w - (w^H V_1[:,i])V_1[:,i]$
	\EndFor
	\If {$||w|| = 0$}
		\State stop
	\Else
		\State $V_1 = [V_1; w/||w||_2]$
	\EndIf
\EndFor
\State $V_{\tilde{\mathcal{G}}_1} = [V_1[:,1:m], V_2[:,2:m]]$
\For {$j = m+1$ to $2m-1$}
	\For {$i = 1$ to $j$}
	\State $V_{\tilde{\mathcal{G}}_1}[:, j] = V_{\tilde{\mathcal{G}}_1}[:, j] - (V_{\tilde{\mathcal{G}}_1}[:, j]^H V_{\tilde{\mathcal{G}}_1}[:, i])V_{\tilde{\mathcal{G}}_1}[:, i]$
	\EndFor
\EndFor
\end{algorithmic}
\caption{The TGSAR I procedure} 
\label{algorithm:Arnoldi3}
\end{algorithm}

\begin{algorithm}[H]
\begin{algorithmic}[1] 
\State $v_1 = b/||b||_2$
\State $w = Av_1$
\State $w = w - (w^H v_1)v_1$
\State $v_2 = w/||w||_2$
\State $V = [v_1, v_2]$
\For {$j = 2$ to $2m-1$}
	\If {$mod(j,2)=0$}, where $mod(a,2)$ is used to find the remainder number after dividing $a$ by $2$.
		\State $w=B V[:,j-2]$
	\Else
		\State $w=A V[:,j-2]$
	\EndIf
	\State $\alpha = ||w||_2$
	\For {$i = 1$ to $j$}
		\State $w = w - (w^H V[:,i])V[:,i]$
	\EndFor
	\If {$||w||_2 < \eta*\alpha$}
		\For {$i = 1$ to $j$}
			\State $w = w - (w^H V[:,i])V[:,i]$
		\EndFor
	\EndIf
	\If {$||w||_2 = 0$}
		\State stop
	\Else
		\State $V = [V; w/||w||_2]$
	\EndIf
\EndFor
\State $V_{\tilde{\mathcal{G}}_1} = V[:,1:2m-1]$
\end{algorithmic}
\caption{The TGSAR II procedure with the partial reorthogonalization technique} 
\label{algorithm:Arnoldi4}
\end{algorithm}

\subsection{The Existing Methods}\label{section:sec3.4}
The reduction subspace of the Arnoldi method is:
$$\mathcal{S}_m(\begin{bmatrix} 0 & I \\ B & A \end{bmatrix};b) = span\{b, \begin{bmatrix} 0 & I \\ B & A \end{bmatrix}b, ..., \begin{bmatrix} 0 & I \\ B & A \end{bmatrix}^{m-1}b\} = span\{U_0 \}$$
where the orthonormal basis $U_0$ is generated by the Arnoldi procedure, and then we can reduce the inverse matrix $\begin{bmatrix} 0 & I \\ B & A \end{bmatrix}$ as $U_0^H \begin{bmatrix} 0 & I \\ B & A \end{bmatrix} U_0$ by imposing the Galerkin condition.

The reduction subspace of the SOAR method is $\mathcal{Q}^1_m$ with $r_{-1} = 0$, and its orthonormal basis $U_1$ is generated by the SOAR procedure; the reduction subspace of the GSOAR method is $\mathcal{Q}^1_m$, and its orthonormal basis $U_2$ is generated by the SOAR procedure; the reduction subspace of the TOAR method is $\mathcal{Q}^1_m$, and its orthonormal basis $U_3$ is generated by the TOAR procedure.
After we obtain the orthonormal basis $U_i$, where $i=1,2,3$, we can reduce the original matrices $A$ and $B$ as $U_i^HAU_i$ and $U_i^HBU_i$ by imposing the Galerkin condition. 

\section{Numerical Examples}\label{section:sec4}
In this section, we compare the LQAR in Algorithm~\ref{algorithm:Arnoldi1}, QAR in  Algorithm~\ref{algorithm:Arnoldi2}, TGSAR in Algorithm~\ref{algorithm:Arnoldi4}, Arnoldi \cite{bai2004second}, SOAR \cite{bai2004second}, and TOAR  procedures \cite{lu2016stability}. Also, six Arnoldi-type methods all apply the partial reorthogonalization technique with $\eta = \sqrt{2}/2$. 
For solving the linearized QEP, we call the QR function in the Python NumPy library, which is based on the BLAS (Basic Linear Algebra Subprograms) library. 
The experiments performed on Win10 with 64-bit operating system and Intel(R) Core(TM) i5-6300HQ CPU @ 2.30GHz 2.30 GHz and RAM 12.0 GB and Python 3.9.2.

\subsection{Examples Description}
First, we discretize a turbofan engine LP rotor into the FE matrices or FE model by the FE procedure \cite{friswellurl}, where the LP shaft is supported by four rolling-element bearings. The FE model has 199 nodes, where each node has 4 degrees of freedom, hence the dimension of all FE matrices is 796, or we can say that the degree of freedom of the FE model is 796. 
Next, we use the FE matrices to construct two examples, which the first is related to modal analysis and the second is related to critical speed analysis. 

\begin{table}[ht]
\caption{The matrices $M$, $C$, $K$ of the two examples}
\centering
\begin{tabular}{cl}
\hline
\textbf{examples} & \textbf{description} \\ \hline
\textbf{1} &
  \begin{tabular}[c]{@{}l@{}}$M=M_0$, $C=C_0+\Omega C_1$, $K=K_0$, $C_0 = \alpha M_r + \beta K_r$, \\ $\Omega=500$, $\alpha=10$, $\beta=1e-5$\end{tabular} \\
\textbf{2} &
  \begin{tabular}[c]{@{}l@{}}$\hat{M}=-n^2 M_0 + j n C_1$, $\hat{C}=jnC_0$, $\hat{K}=K_0$,\\ $C_0 = \alpha M_r + \beta K_r$, $n=1$, $\alpha=10$, $\beta=1e-5$\end{tabular} \\ \hline
\end{tabular}
\label{tab:Examples}
\end{table}

\begin{table}[ht]
\caption{The properties of the examples' $M$, $C$, $K$}
\centering
{
\begin{tabular}{clclclc}
\hline
\textbf{LP Rotor} & \multicolumn{2}{c}{\textbf{$M$}}    & \multicolumn{2}{c}{\textbf{$C$}}     & \multicolumn{2}{c}{\textbf{$K$}} \\ \cline{2-7} 
$796 \times 796$ & \multicolumn{1}{c}{\textbf{sparsity}} & \textbf{symmetry} & \multicolumn{1}{c}{\textbf{sparsity}} & \textbf{symmetry} & \multicolumn{1}{c}{\textbf{sparsity}} & \textbf{symmetry} \\ \hline
\textbf{example1}         & 0.992532 & symmetry     & 0.985064 & asymmetry & 0.992538  & symmetry \\
\textbf{example2}          & 0.985064 & asymmetry & 0.992532 & symmetry & 0.992538  & symmetry \\ \hline
\end{tabular}
}
\label{tab:LPMCK}
\end{table}

\begin{figure}[H]
\centering
\includegraphics[width=0.9\textwidth]{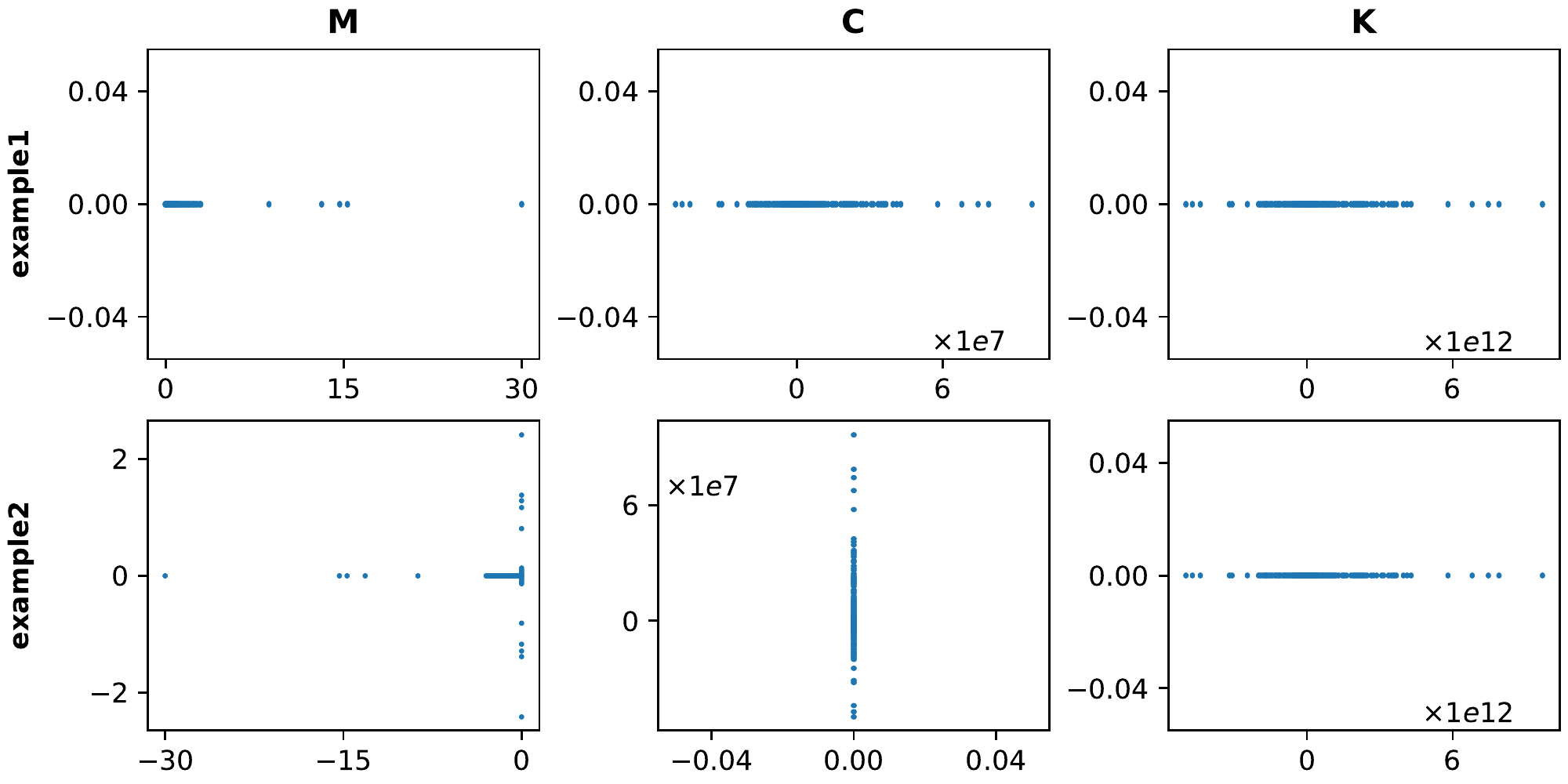}
\caption{The elements of the examples' M, C, K in the complex plane} 
\label{fig:LPMCK}
\end{figure}

Table~\ref{tab:Examples} shows the description of two examples, where $M_r$ is the rotor shaft mass matrix, $K_r$ is the rotor shaft stiffness matrix, $C_b$ is the bearing damping matrix, $\alpha$ and $\beta$ are the Rayleigh damping cofficients. The properties of $M,C,K$ are shown in Table~\ref{tab:LPMCK}. Moreover, Fig.~\ref{fig:LPMCK} shows the elements distribution of $M,C,K$ in the complex plane.

\subsection{Comparisons}
We set the start vectors $b = b_0 = b_1 = 1$, and then compare first $h$ approximate conjugate eigenvalues of two examples computed by six reduction methods. It also can be considered that we reduce the FE model of 796 degrees of freedom to $\ceil*{h/2}$ degrees of freedom. 

\begin{table}[ht]
\caption{The first 10 eigenvalues computed by six methods with $m=10$}
\centering
\resizebox{\textwidth}{!}{
\begin{tabular}{|c|c|c|c|c|c|c|c|}
\hline
\textbf{example3} & \textbf{Exact} & \textbf{Arnoldi} & \textbf{SOAR} & \textbf{TOAR} & \textbf{LQAR} & \textbf{QAR}  & \textbf{TGSAR} \\ \hline
\textbf{1}        & -3.03-279.08j  & -3.22-279.16j    & -3.03-279.08j & -2.99+279.16j & -3.05-279.11j & -3.03+279.07j & -3.03-279.08j  \\ \hline
\textbf{2}        & -3.03+279.08j  & -3.22+279.16j    & -3.03+279.08j & -2.99-279.16j & -3.05+279.11j & -3.03-279.07j & -3.03+279.08j  \\ \hline
\textbf{3}        & -4.17-392.98j  & -3.08+393.14j    & -4.19+392.97j & -4.18-393.12j & -4.12+392.9j  & -4.22-392.99j & -4.17-392.96j  \\ \hline
\textbf{4}        & -4.17+392.98j  & -3.08-393.14j    & -4.19-392.97j & -4.18+393.12j & -4.12-392.9j  & -4.22+392.99j & -4.17+392.96j  \\ \hline
\textbf{5}        & -5.38-458.51j  & -13.78+494.65j   & -4.94+460j    & -5.15-458.45j & -5.36+459.27j & -5.36-458.7j  & -5.42-458.5j   \\ \hline
\textbf{6}        & -5.38+458.51j  & -13.78-494.65j   & -4.94-460j    & -5.15+458.45j & -5.36-459.27j & -5.36+458.7j  & -5.42+458.5j   \\ \hline
\textbf{7}        & -5.16-492.25j  & -267.52+970.07j  & -5.23-495.11j & -5.43+494.26j & -5.39+492.55j & -5.14-492.23j & -5.16+492.25j  \\ \hline
\textbf{8}        & -5.16+492.25j  & -267.52-970.07j  & -5.23+495.11j & -5.43-494.26j & -5.39-492.55j & -5.14+492.23j & -5.16-492.25j  \\ \hline
\textbf{9}        & -5.51-497.63j  & -1148.78-0j      & -6.63-522.52j & -6.84-515.44j & -5.28+497.61j & -5.38-498.04j & -5.51+497.63j  \\ \hline
\textbf{10}       & -5.51+497.63j  & 852-980.69j      & -6.63+522.52j & -6.84+515.44j & -5.28-497.61j & -5.38+498.04j & -5.51-497.63j  \\ \hline
\textbf{error}    & -              & 4542.96          & 61.86         & 43.98         & 3.49          & 1.72          & 0.14           \\ \hline
\end{tabular}
}
\label{tab:compare-example1}
\end{table}

\begin{table}[ht]
\caption{The first 10 eigenvalues computed by six methods with $m=10$}
\centering
\resizebox{\textwidth}{!}{
\begin{tabular}{|c|c|c|c|c|c|c|c|}
\hline
\textbf{example3} & \textbf{Exact} & \textbf{Arnoldi} & \textbf{SOAR}       & \textbf{TOAR}         & \textbf{LQAR}             & \textbf{QAR}              & \textbf{TGSAR}            \\ \hline
\textbf{1}        & -325.21+3.13j  & -0.22-35.9j      & 186.37+0.91j        & -156.31+2.79j         & -325.21+3.13j             & -325.21+3.13j             & -325.21+3.12j             \\ \hline
\textbf{2}        & 325.21+3.13j   & -324.41+1.84j    & -186.74+6.77j       & 156.38+7.97j          & 325.22+3.14j              & 325.21+3.13j              & 325.21+3.14j              \\ \hline
\textbf{3}        & 410.97+3.56j   & 324.45+4.43j     & -359.33+38.94j      & -749.89+27.76j        & 411+3.53j                 & -410.97+3.56j             & 410.95+3.56j              \\ \hline
\textbf{4}        & -410.97+3.56j  & -437.87+2.54j    & 361.42-29.46j       & 751.1-15.28j          & -411+3.59j                & 410.97+3.56j              & -410.95+3.56j             \\ \hline
\textbf{5}        & 455.42+5.85j   & 439.43+3.94j     & 10.73+671.99j       & 1149.9+1768.59j       & -455.16+5.69j             & 455.4+5.91j               & -455.5+5.85j              \\ \hline
\textbf{6}        & -455.42+5.85j  & -1.74-509.61j    & -9.94-701.18j       & -1160.1-1773.75j      & 455.18+6.04j              & -455.41+5.79j             & 455.5+5.85j               \\ \hline
\textbf{7}        & -492.5+4.99j   & 1.86+526.16j     & 707.54+201.83j      & -16369.19+46902.76j   & 492.69+4.97j              & 492.5+5j                  & 492.5+4.98j               \\ \hline
\textbf{8}        & 492.5+4.99j    & -529.03-20.01j   & -718.73-178.02j     & 16427.42-46909.35j    & -492.69+5.02j             & -492.5+4.98j              & -492.5+4.99j              \\ \hline
\textbf{9}        & -497.43+5.93j  & 531.22+30.1j     & 85407.23+11838.87j  & 414617.11-716010.22j  & 498.1+6.16j               & 497.44+5.95j              & 497.56+5.96j              \\ \hline
\textbf{10}       & 497.43+5.93j   & 4.67+1848.51j    & -85412.75-11807.69j & -414678.47+716012.45j & -498.11+5.75j             & -497.44+5.91j             & -497.56+5.92j             \\ \hline
\textbf{error}    & -              & 4907.11          & 196968.21           & 2391912.96            & 3.16 & 0.22 & 0.54 \\ \hline
\end{tabular}
}
\label{tab:compare-example2}
\end{table}

First, we set $m=10$ and compute the first 10 approximation eigenvalues. The reuslts are shown in Table~\ref{tab:compare-example1} and Table~\ref{tab:compare-example2}.
We define an error function as:
\begin{equation}
error = \sum_{i=1}^{10}abs(abs(Re(s_i^e)) - abs(Re(s_i^a))) + \sum_{i=1}^{10}abs(abs(Im(s_i^e)) - abs(Im(s_i^a)))
\end{equation}
where $s^e$ means exact eigenvalue, and $s^a$ means approximation eigenvalue. And, to avoid confusion, we use $abs(\cdot)$ instead of $|\cdot|$ for absolute function.
In example1, we can rank the six methods by errors as: TGSAR $>$ QAR $>$ LQAR $>$ TOAR $>$ SOAR $>$ Arnoldi. 
Similarly, in example2, the ranking is: QAR $>$ TGSAR $>$ LQAR $>$ Arnoldi $>$ SOAR $>$ TOAR. The results in both examples show the superior accuracy of our proposed three Arnoldi-type methods.

\begin{figure}[H]
\centering
\includegraphics[width=0.9\textwidth]{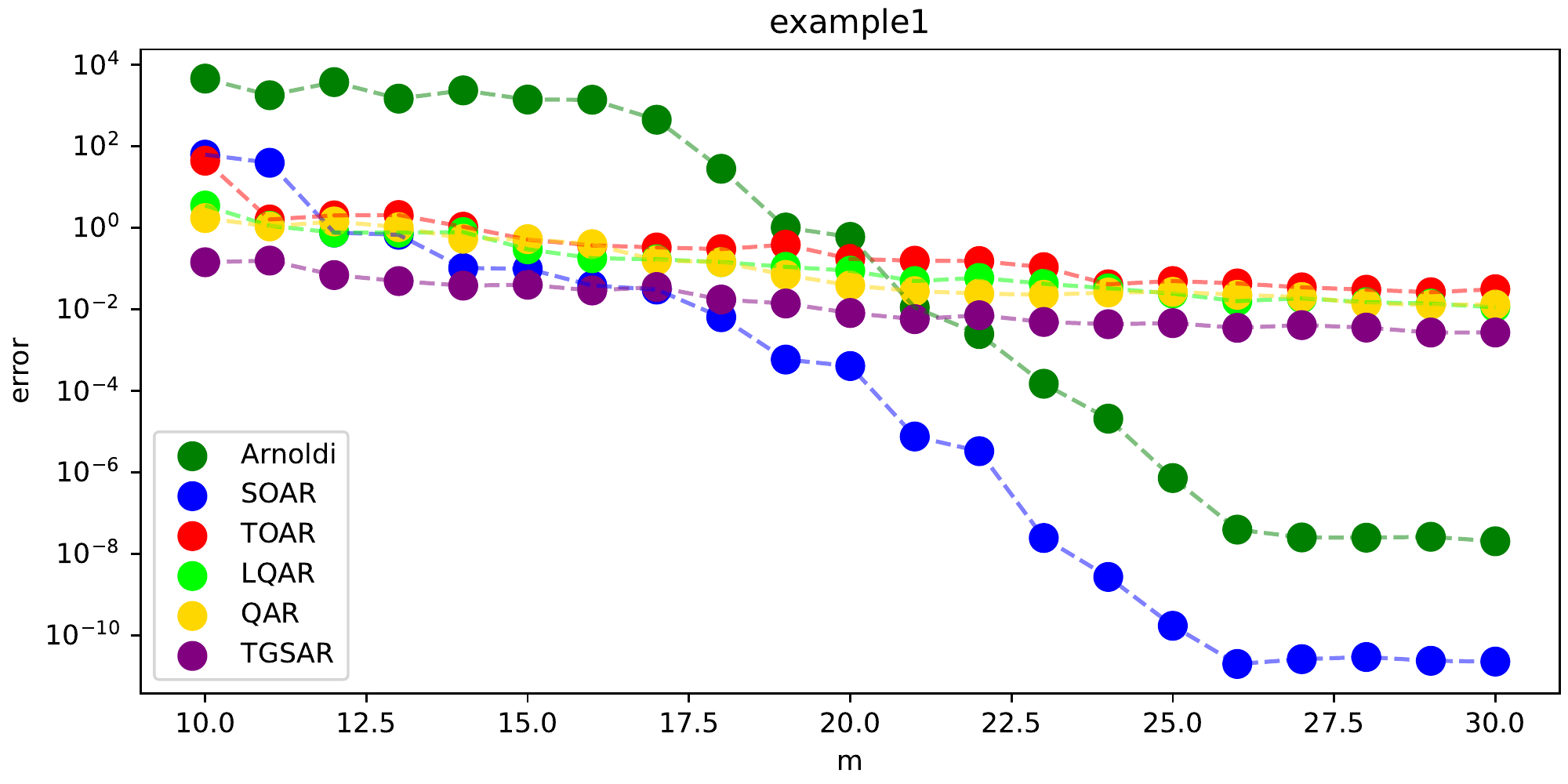}
\caption{The comparison of computational accuracy of the first 10 eigenvalues computed by six methods} 
\label{fig:example11}
\end{figure}

\begin{figure}[H]
\centering
\includegraphics[width=0.9\textwidth]{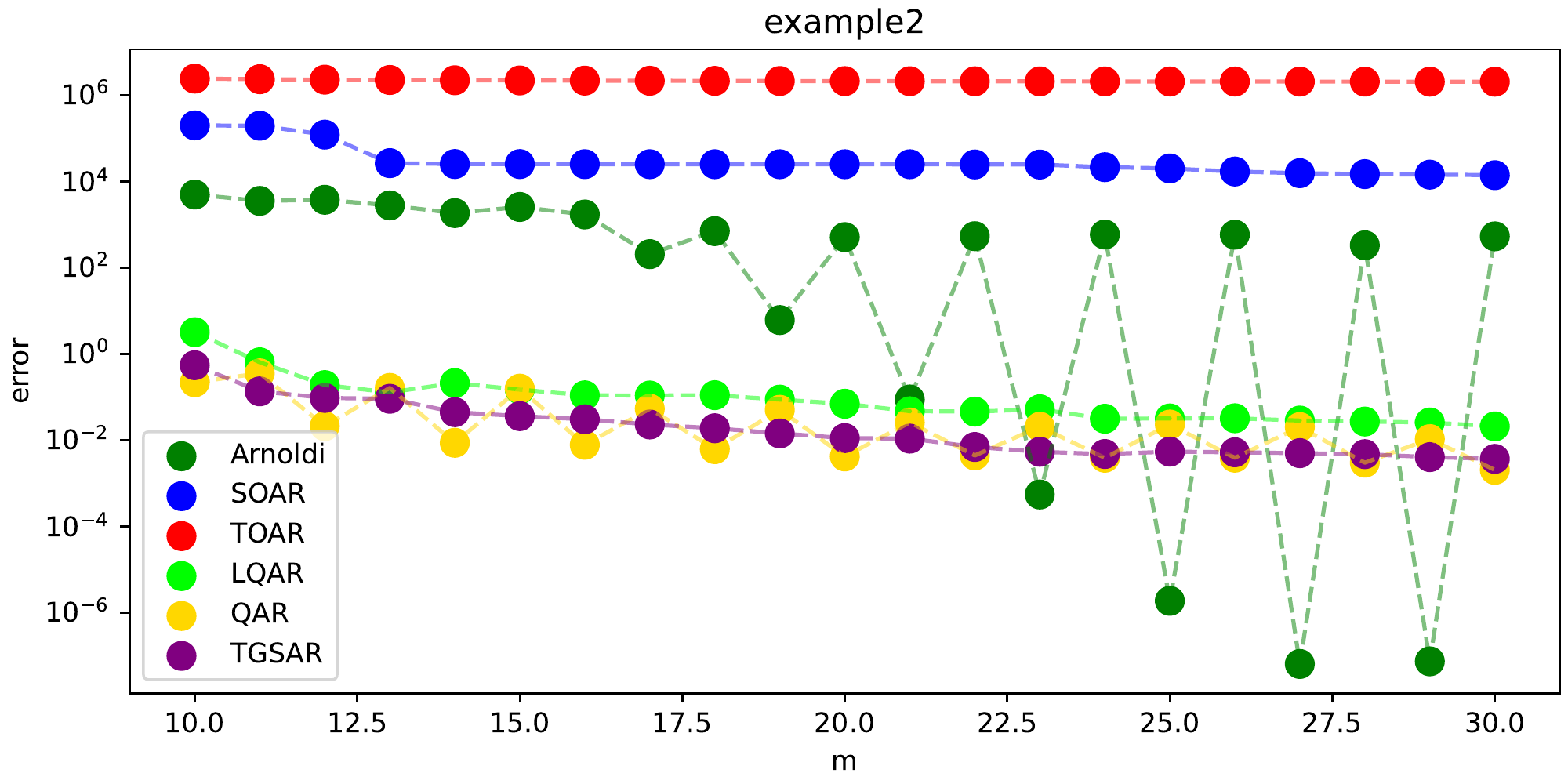}
\caption{The comparison of computational accuracy of the first 10 eigenvalues computed by six methods} 
\label{fig:example12}
\end{figure}

Next, we increase $m$ and compare the errors of the first 10 approximation eigenvalues by six methods. The results are shown in Fig~\ref{fig:example11} and Fig~\ref{fig:example12}. 
In example1, we can see that when m increases, the Arnoldi, SOAR, and TOAR will improve their accuracy, and the first two methods eventually outperform than our proposed three methods.
In example2, we can see that the Arnoldi is sometimes better than our proposed three methods as m increases, but it is unstable.

\begin{figure}[ht]
\centering
\includegraphics[width=0.9\textwidth]{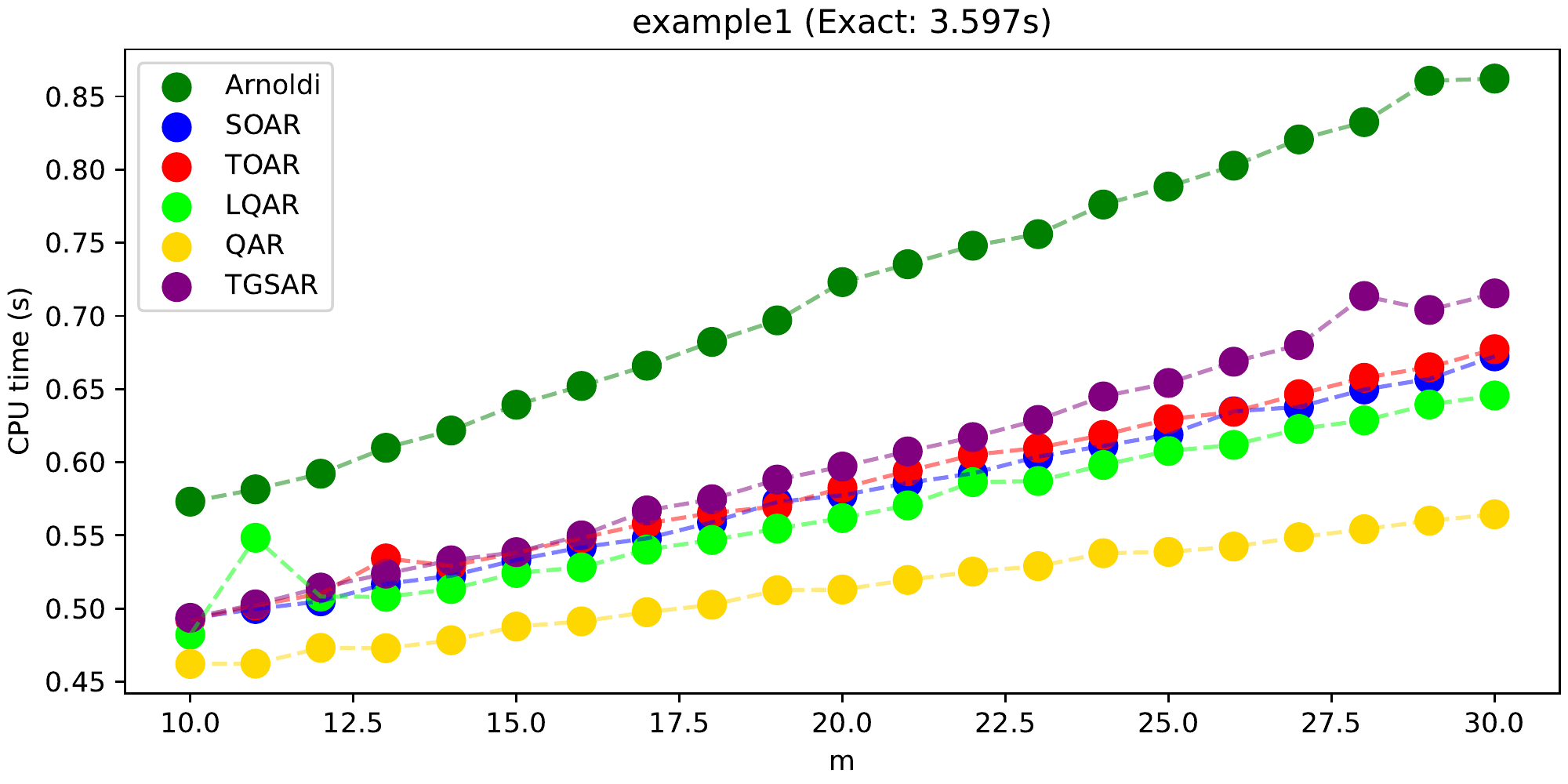}
\caption{The comparison of computational effiency of the first 10 eigenvalues computed by six methods} 
\label{fig:example21}
\end{figure}

\begin{figure}[ht]
\centering
\includegraphics[width=0.9\textwidth]{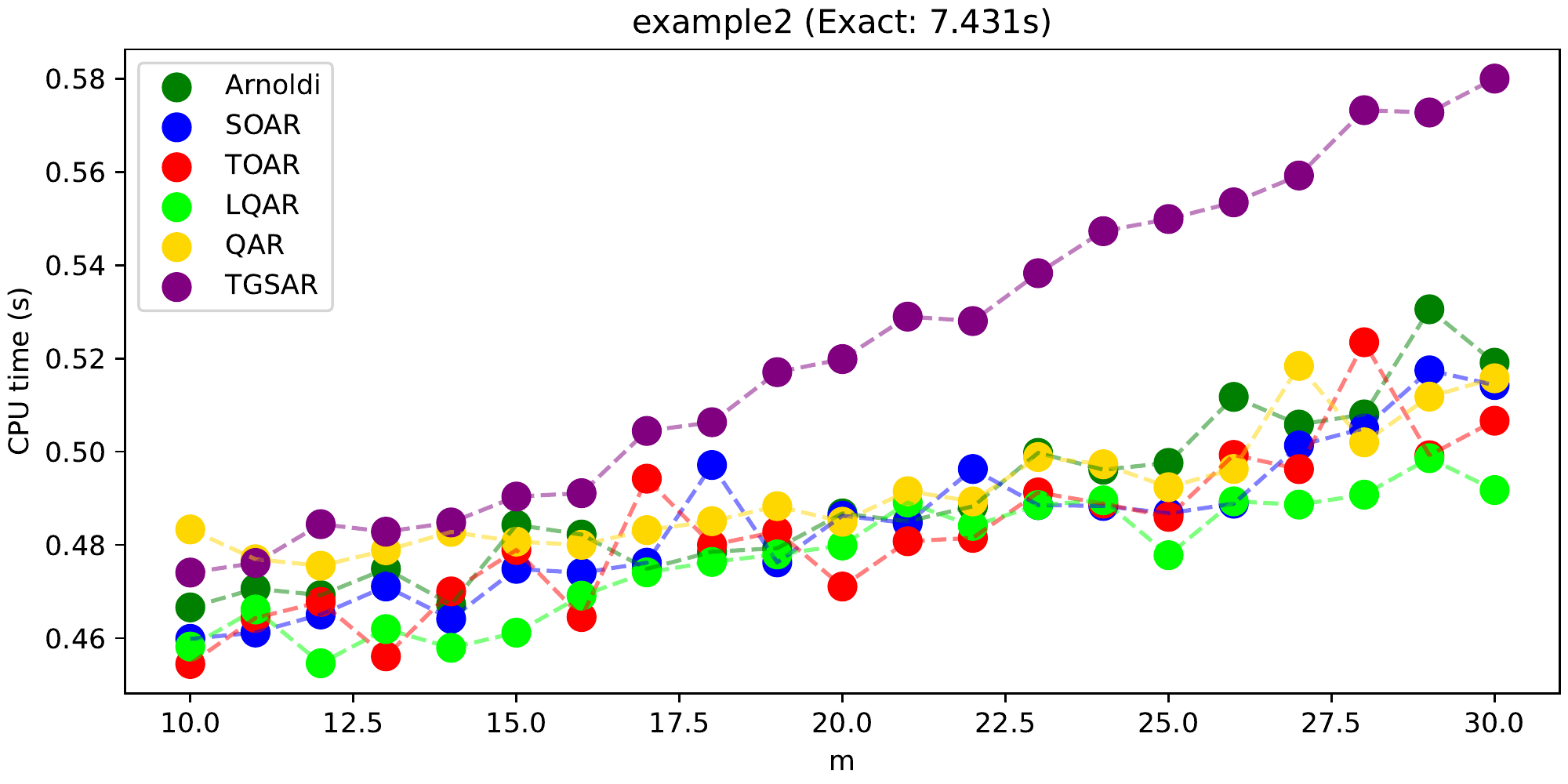}
\caption{The comparison of computational effiency of the first 10 eigenvalues computed by six methods} 
\label{fig:example22}
\end{figure}

After that, we compare the CPU time required by all methods. The results are shown in Fig~\ref{fig:example21} and Fig~\ref{fig:example22}, where each data is the average of running 10 times. As can be seen from Fig~\ref{fig:example21}, the six methods decrease the 3.597s required for the exact solution to less than 0.9s. 
Similarly, Fig~\ref{fig:example22} shows that the six methods reduce the 7.431s required for the exact solution to less than 0.6s. 
In addition, the Arnoldi method in Fig~\ref{fig:example21} is more time-consuming than the other methods, because although Arnoldi needs to perform a fewer number of the modified Gram-Schmidt iterations, it takes more time to do Ritz transformation.

\begin{figure}[H]
\centering
\includegraphics[width=0.9\textwidth]{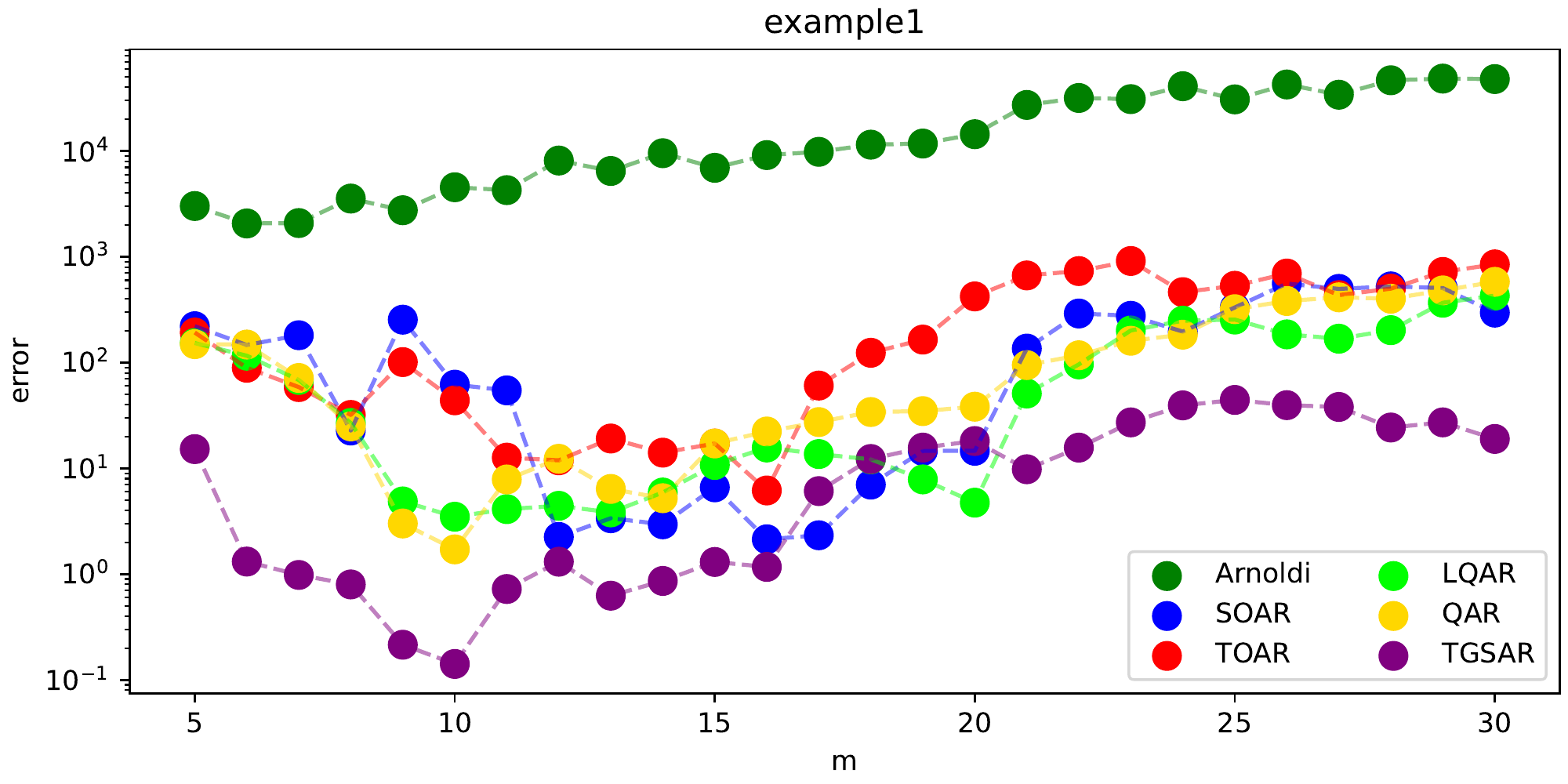}
\caption{The comparison of computational accuracy of the first $m$ eigenvalues computed by six methods} 
\label{fig:example31}
\end{figure}

\begin{figure}[H]
\centering
\includegraphics[width=0.9\textwidth]{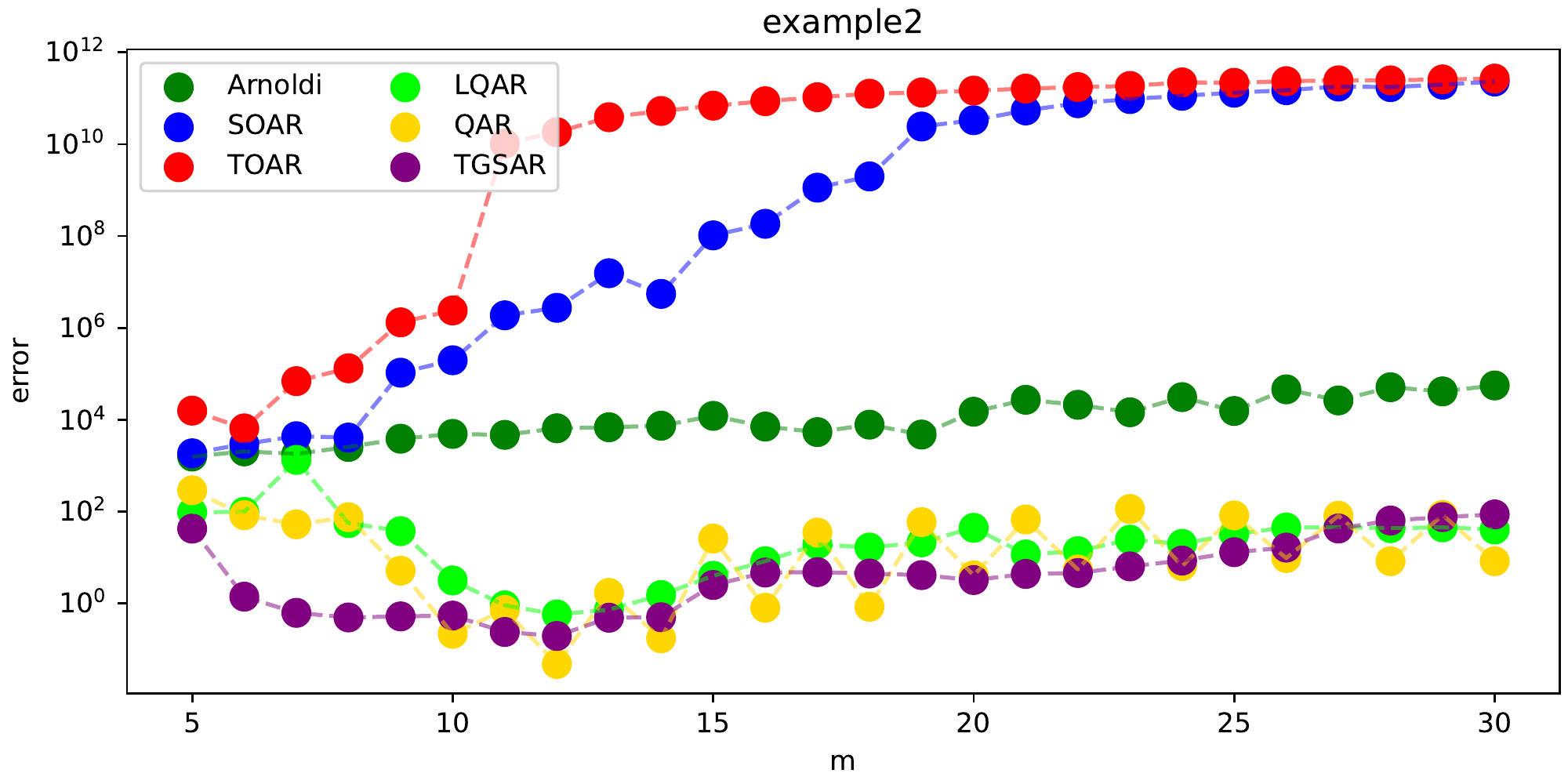}
\caption{The comparison of computational accuracy of the first $m$ eigenvalues computed by six methods} 
\label{fig:example32}
\end{figure}

Finally, we increase $m$, and compre the errors of the first $m$ approximation eigenvalues. The results are shown in Fig~\ref{fig:example31} and Fig~\ref{fig:example32}. 
We define a score function as:
\begin{equation}
score = \sum_{i}^5 n_i*(1-(i-1)*0.2)
\end{equation}
where $n_i$ is the number of times when one method ranks as $i$th in $m$ error comparisons.
In example1, these methods can be ranked by scores as: TGSAR $>$ LQAR $>$ SOAR $>$ QAR $>$ TOAR $>$ Arnoldi. In example2, the ranking is: TGSAR $>$ QAR $>$ LQAR $>$ Arnoldi $>$ SOAR $>$ TOAR. As a result, from comprehensive consideration of computation efficiency and accuracy, our proposed three methods are superior to the existing methods.

\subsection{Discussions}
The linearized QEP Krylov subspace $\mathcal{Q}_m^1$ and the SEP Krylov subspace $\mathcal{S}_m(\begin{bmatrix} 0 & I \\ B & A \end{bmatrix};b)$ have the relationship as:
\begin{equation}
\mathcal{S}_m(\begin{bmatrix} 0 & I \\ B & A \end{bmatrix};b) \subseteq span\{\begin{bmatrix} V_{\mathcal{Q}^1} & 0 \\ 0 & V_{\mathcal{Q}^1} \end{bmatrix} \}
\end{equation}
Therefore, the LQAR, SOAR, and TOAR methods which generate the orthonormal basis of the subspace $\mathcal{Q}_m^1$ are at least as good as the Arnoldi method which generates the orthonormal basis of the subspace $\mathcal{S}_m(\begin{bmatrix} 0 & I \\ B & A \end{bmatrix};b)$. 
However, in the numerical examples, only the LQAR method is always better than the Arnoldi method. 
This means the Arnoldi iterative formulation is better than the iterative formulations of the SOAR and TOAR, and can generate the more accurate orthonormal basis of $\mathcal{Q}_m^1$. It is unclear to us where the problems with the iterative formulations of the SOAR and TOAR are and requires further study.  In addition, due to our proposed three new methods all use the Arnoldi iteration formulation, and then they can generate the more accurate orthonormal basis of $\mathcal{Q}_m^1$, $\mathcal{Q}_m^2$ and $\tilde{\mathcal{G}}_{2m-1}$.

The dimension of the generalized SEP Krylov subspace $\mathcal{G}_{2^m-1}$ is large, and then we present that it only need to be truncated as $\tilde{\mathcal{G}}_{2m-1}$, and now we explain the reason.
In section~\ref{section:sec3.1}, we know that the essence of reducing the QEP is to find a subspace $\mathcal{Q}_m$ that can reduce the matrices $A$ and $B$ at the same time. And, the SEP Krylov subspaces $\mathcal{S}_m(A; b)$ and $\mathcal{S}_m(B; b)$ can be used to reduce the matrices $A$ and $B$, respectively. Therefore, at least we should truncate $\mathcal{G}_{2^m-1}$ as:
$$\tilde{\mathcal{G}}_{2m-1} = \mathcal{S}_m(A; b) \cup \mathcal{S}_m(B; b)$$
In practice, $\tilde{\mathcal{G}}_{2m-1}$ is enough for reducing the QEP in rotor dynamics, and we give another numerical examples in \ref{section:another}.
Since $\mathcal{Q}^2_m$ have a relationship with $\mathcal{G}_{2^m-1}$:
$$\mathcal{Q}^2_m \subseteq \mathcal{G}_{2^m-1}$$
Therefore, $\tilde{\mathcal{G}}_{2m-1}$ is an approximation to the space including $\mathcal{Q}_m^2$.

In this paper, all reduction subspaces $\mathcal{S}_m(\begin{bmatrix} 0 & I \\ B & A \end{bmatrix};b)$, $\mathcal{Q}_m^1$, $\mathcal{Q}_m^2$, and $\tilde{\mathcal{G}}_{2m-1}$ have the relationship as:
\begin{enumerate}[1)]
	\item $\mathcal{S}_m(\begin{bmatrix} 0 & I \\ B & A \end{bmatrix};b)$ is an approximation to $span\{\begin{bmatrix} V_{\ceil*{m/2}} & 0 \\ 0 & V_{\floor*{m/2}} \end{bmatrix} \}$
	\item $\mathcal{Q}^1_m, \mathcal{Q}^2_m$, and $\tilde{\mathcal{G}}_{2m-1}$ are approximations to $span\{V_m\}$
	\item $\tilde{\mathcal{G}}_{2m-1}$ is an approximation to the space including $span\{V_m\}$
	\item $\mathcal{S}_m(\begin{bmatrix} 0 & I \\ B & A \end{bmatrix};b) \subseteq span\{\begin{bmatrix} V_{\mathcal{Q}^1} & 0 \\ 0 & V_{\mathcal{Q}^1} \end{bmatrix} \}$
\end{enumerate}
where $V_m$ is a matrix constructed by first $m$ exact eigenvectors of the original QEP, and $\floor*{\cdot}$ is floor function. 

\section{Concluding remarks}\label{section:sec5}
In this paper, we present three QEP reduction subspaces. 
First is the linearized QEP Krylov subspace $\mathcal{Q}_m^1$, which is spanned by the power iteration sequences of the linearized QEP and already used in the SOAR, GSAOR, and TOAR methods. 
Second is the QEP Krylov subspace $\mathcal{Q}_m^2$, which is spanned by the power iteration sequences of the QEP. 
The last one is the truncated generalized SEP Krylov subspace $\tilde{\mathcal{G}}_{2m-1}$, which is an extension of the second subspace.
Finally, we give three Arnoldi-type procedures to generate the orthonormal basis of these three subspaces.
Numerical examples show the high accuracy of our proposed methods.

The LQAR, QAR, and TGSAR methods we proposed are extensions of the Arnoldi method. These methods can be well applied to accelerate the damped modal analysis and critical speed analysis or to reduce the QEP in rotor dynamics. 
It remains to be determined whether these methods can be extended to reduce the polynomial eigenvalue problem and whether they can be applied to accelerate fluid mechanics analysis, signal processing, and electromagnetics analysis.

\section*{Declaration of competing interest}
The authors declare that they have no known competing financial interests or personal relationships that can have appeared to influence the work reported in this paper.

\section*{Acknowledgements}
The authors would like to thank M. I. Friswell, J. E. Penny, S. D. Garvey, and A. W. Lees for their rotor dynamics FE procedure.

\begin{appendix}

\section{Another Numerical Examples}
\label{section:another}
We use the other two rotors to construct the examples. First is fan rotor which is the first half of the LP rotor, and it is supported by two rolling-element bearings. Second is a small-scale rotor of the rotor test bench, and it is also supported by two rolling-element bearings. Each rotor construct two examples similar to the LP rotor.

\begin{figure}[H]
\centering
\includegraphics[width=0.9\textwidth]{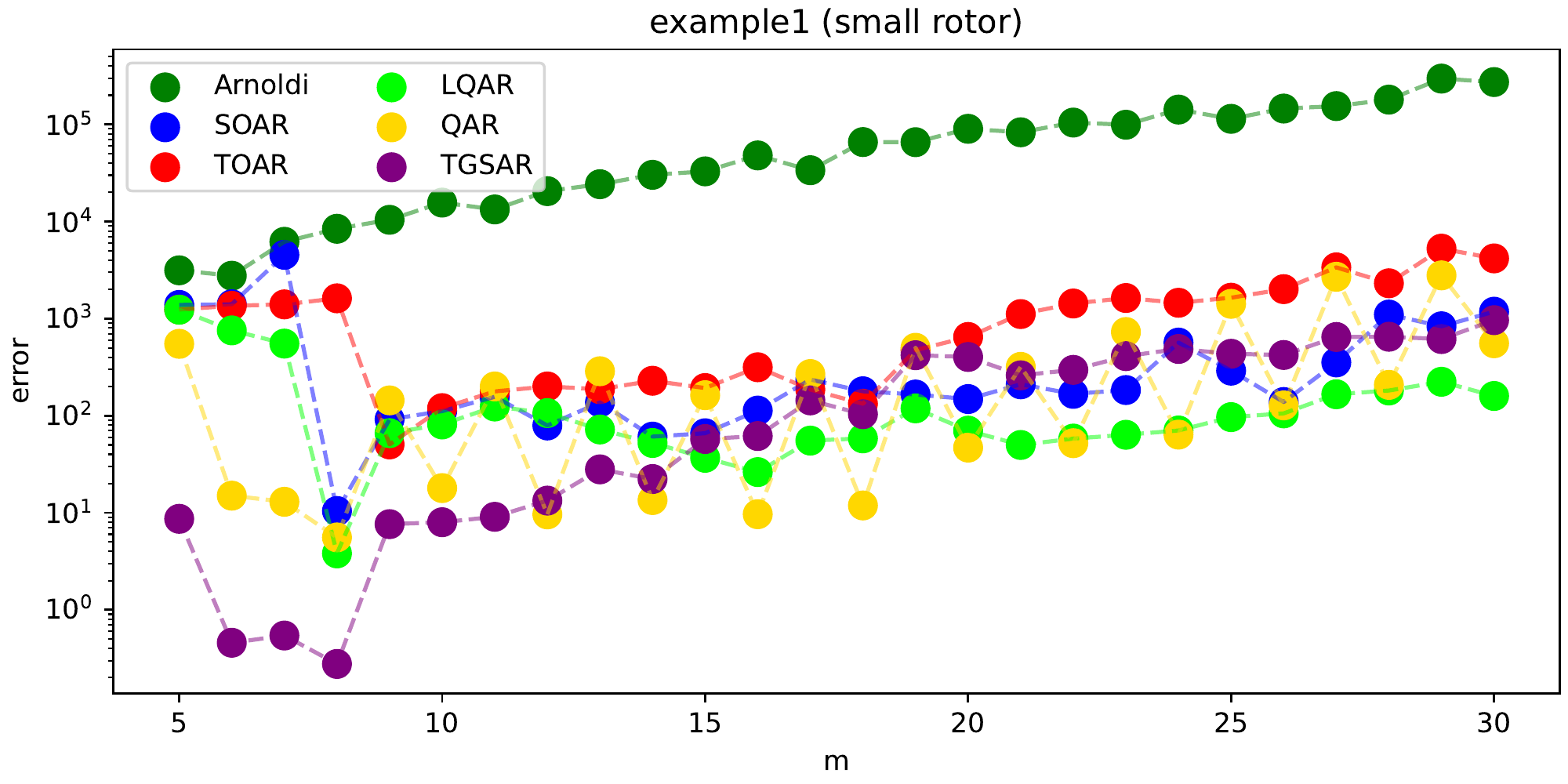}
\caption{The comparison of computational accuracy of the first $m$ eigenvalues computed by six methods} 
\label{fig:small1}
\end{figure}

\begin{figure}[H]
\centering
\includegraphics[width=0.9\textwidth]{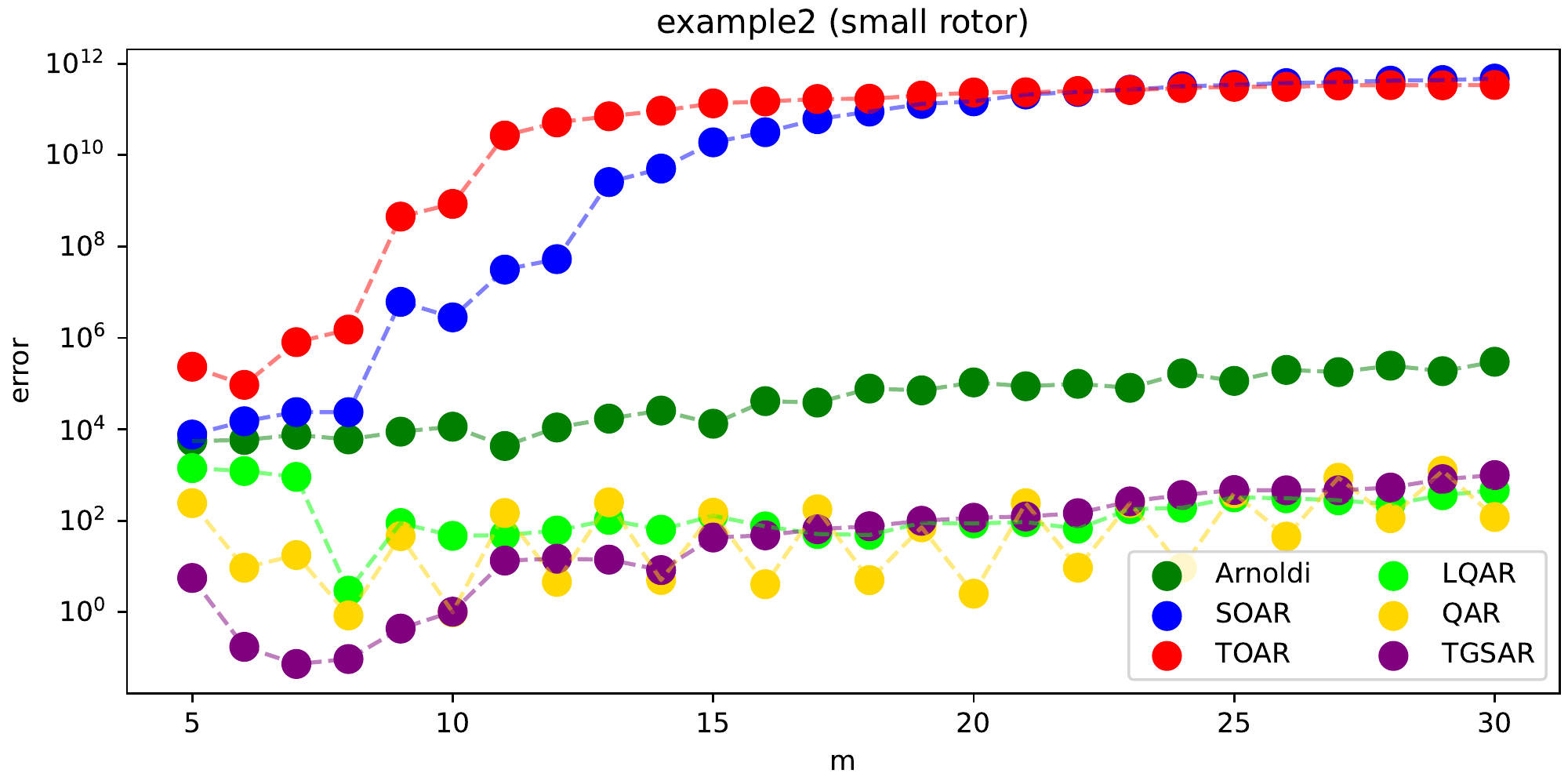}
\caption{The comparison of computational accuracy of the first $m$ eigenvalues computed by six methods} 
\label{fig:small2}
\end{figure}

\begin{figure}[H]
\centering
\includegraphics[width=0.9\textwidth]{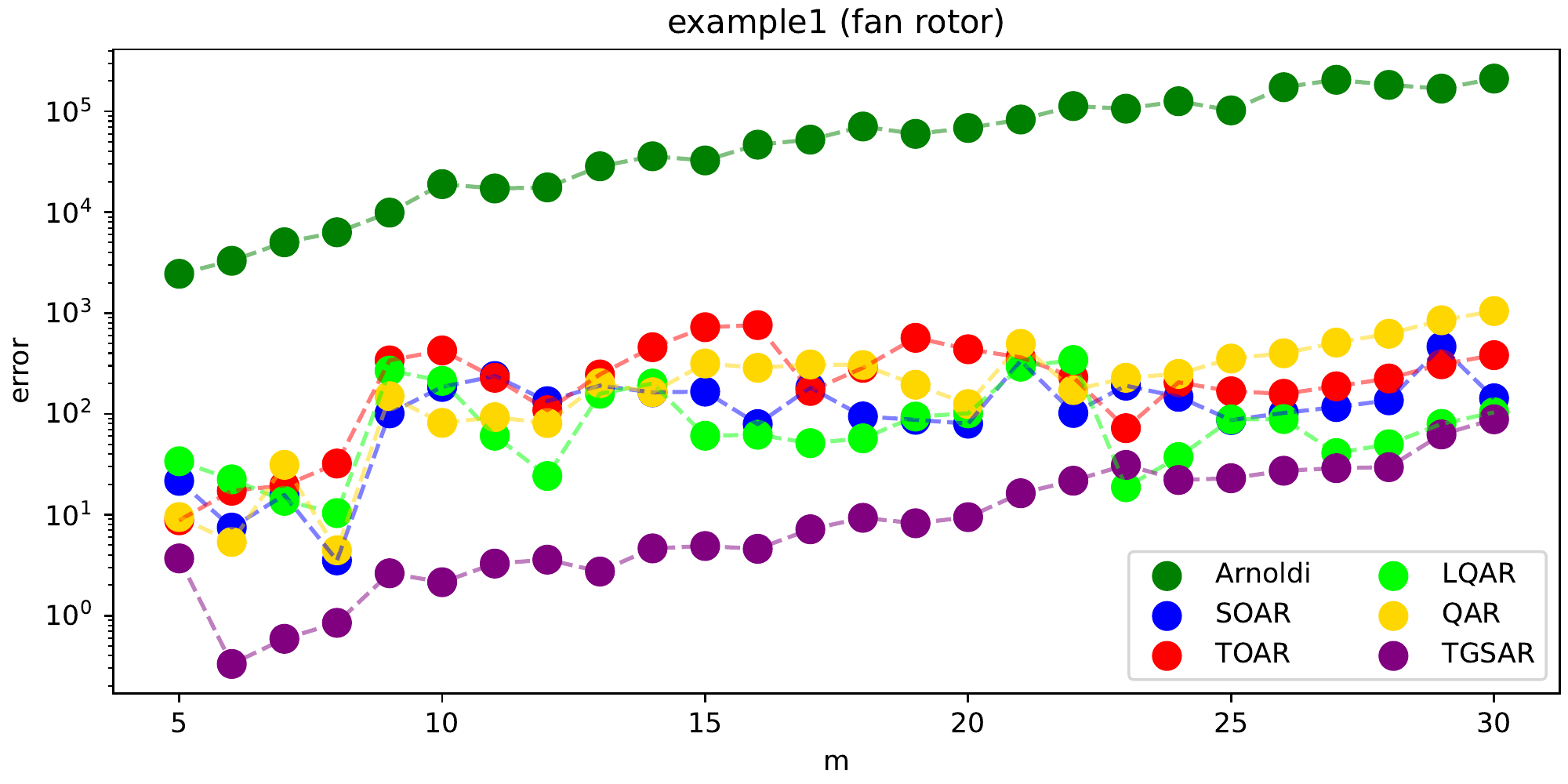}
\caption{The comparison of computational accuracy of the first $m$ eigenvalues computed by six methods} 
\label{fig:fan1}
\end{figure}

\begin{figure}[H]
\centering
\includegraphics[width=0.9\textwidth]{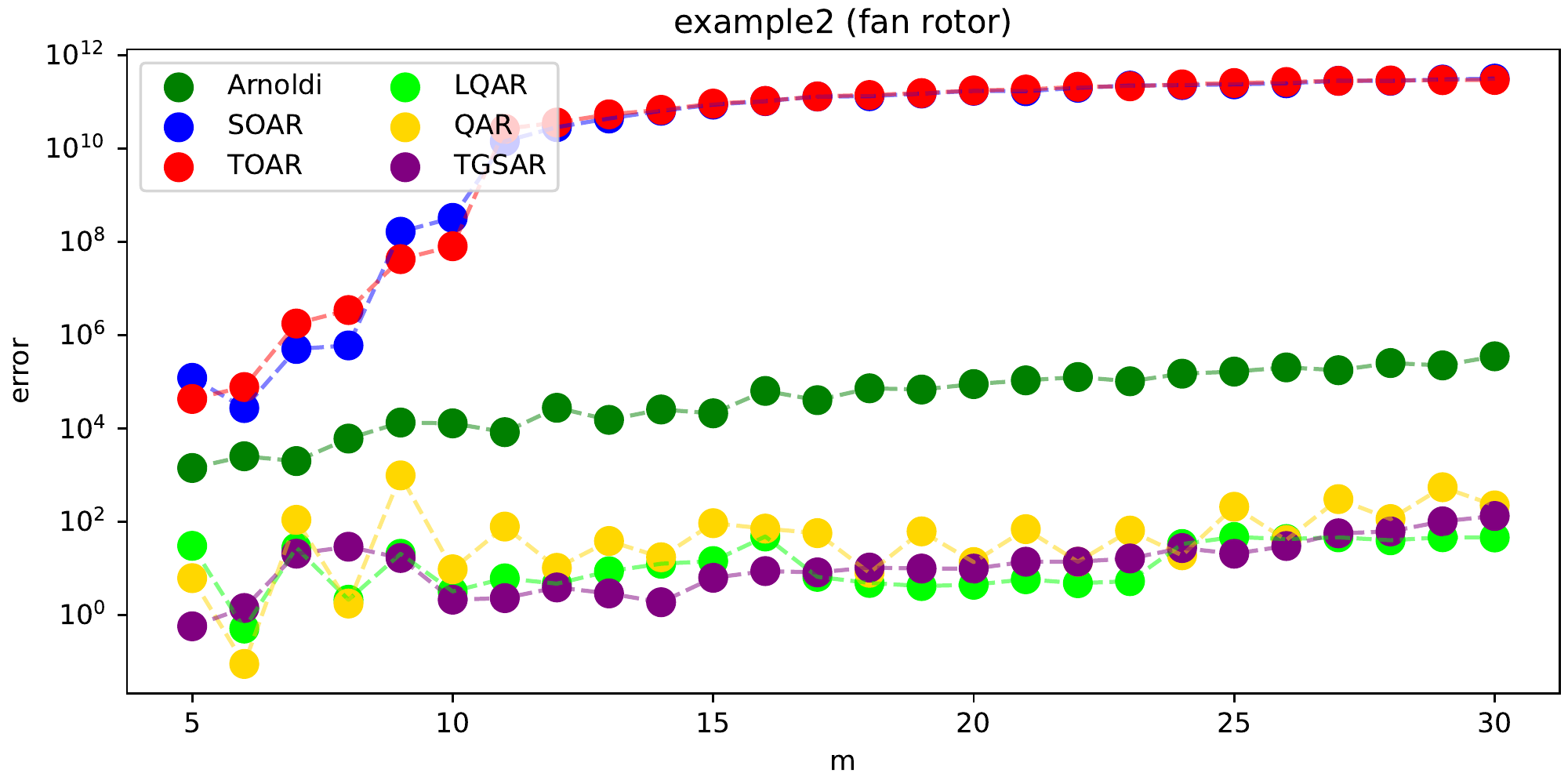}
\caption{The comparison of computational accuracy of the first $m$ eigenvalues computed by six methods} 
\label{fig:fan2}
\end{figure}

We increase $m$, and compre the errors of the first $m$ approximation eigenvalues. The results shown in Fig.~\ref{fig:small1}, Fig.~\ref{fig:small2}, Fig.~\ref{fig:fan1}, and Fig.~\ref{fig:fan2}. We can see that three new methods are still better than the existing methods.

\begin{table}[ht]
\caption{The properties of the examples' $M$, $C$, $K$}
\centering
{
\begin{tabular}{ccccccc}
\hline
\multirow{2}{*}{\textbf{\begin{tabular}[c]{@{}c@{}}Small Rotor\\ $236\times 236$\end{tabular}}} & \multicolumn{2}{c}{\textbf{M}}        & \multicolumn{2}{c}{\textbf{C}}        & \multicolumn{2}{c}{\textbf{K}}        \\ \cline{2-7} 
                                                                                                & \textbf{sparsity} & \textbf{symmetry} & \textbf{sparsity} & \textbf{symmetry} & \textbf{sparsity} & \textbf{symmetry} \\ \hline
\textbf{example1}                                                                               & 0.977521          & symmetry          & 0.955042          & asymmetry         & 0.977593          & symmetry          \\
\textbf{example2}                                                                               & 0.955042          & asymmetry         & 0.977521          & symmetry          & 0.977593          & symmetry          \\ \hline
\multicolumn{7}{l}{Note: $\Omega =200, \beta=10, \alpha=1e-5$}                                                                                                                                                                                              
\end{tabular}
}
\end{table}

\begin{table}[ht]
\caption{The properties of the examples' $M$, $C$, $K$}
\centering
{
\begin{tabular}{ccccccc}
\hline
\multirow{2}{*}{\textbf{\begin{tabular}[c]{@{}c@{}}Fan Rotor\\ $440\times 440$\end{tabular}}} & \multicolumn{2}{c}{\textbf{M}}        & \multicolumn{2}{c}{\textbf{C}}        & \multicolumn{2}{c}{\textbf{K}}        \\ \cline{2-7} 
                                                                                              & \textbf{sparsity} & \textbf{symmetry} & \textbf{sparsity} & \textbf{symmetry} & \textbf{sparsity} & \textbf{symmetry} \\ \hline
\textbf{example1}                                                                             & 0.986446          & symmetry          & 0.972893          & asymmetry         & 0.986446          & symmetry          \\
\textbf{example2}                                                                             & 0.972893          & asymmetry         & 0.986446          & symmetry          & 0.986446          & symmetry          \\ \hline
\multicolumn{7}{l}{Note: $\Omega =400, \beta=10, \alpha=1e-5$}                                                                                                                                                                                              
\end{tabular}
}
\end{table}

\begin{figure}[H]
\centering
\includegraphics[width=0.9\textwidth]{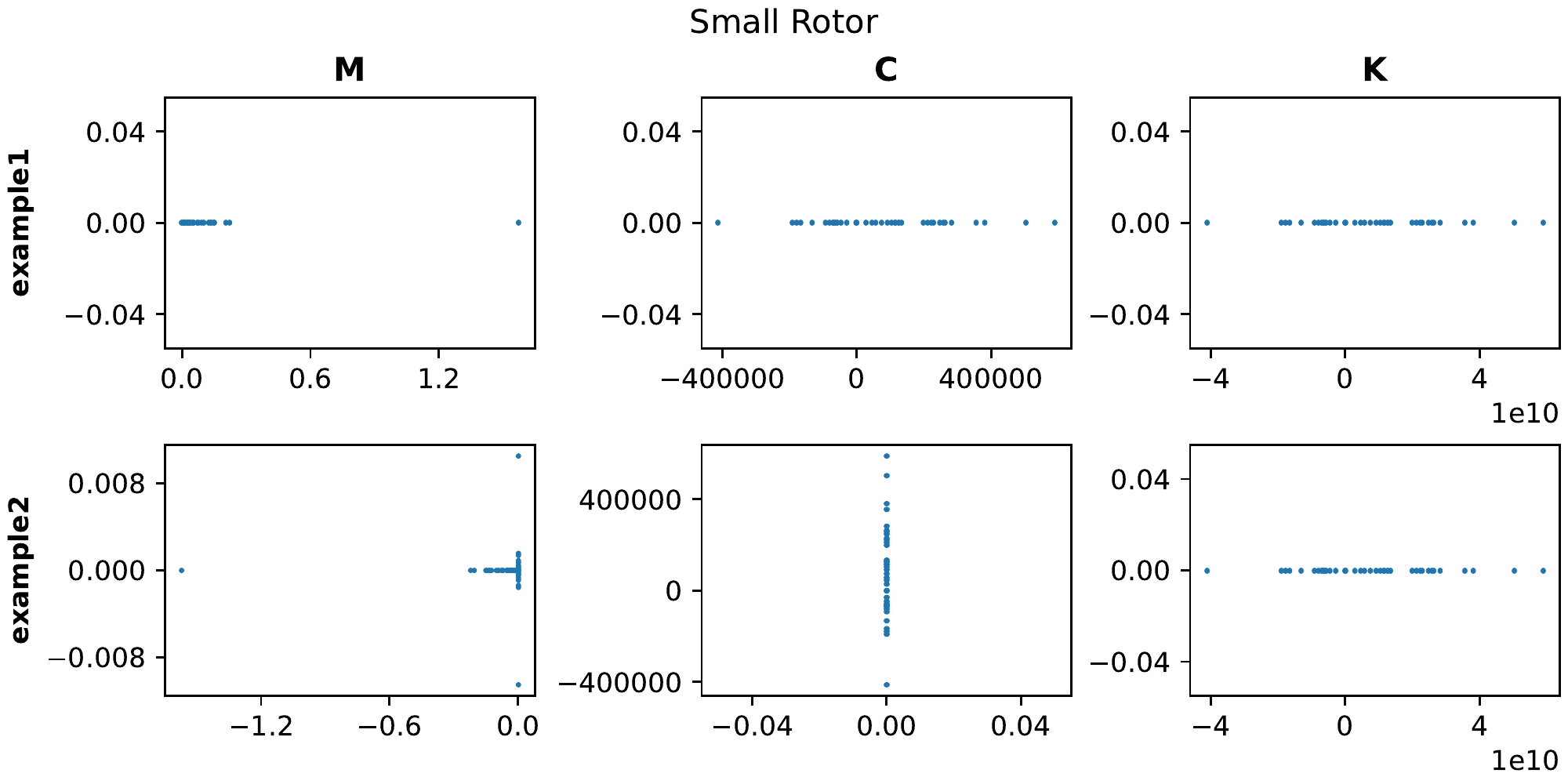}
\caption{The elements of the examples' M, C, K in the complex plane} 
\end{figure}

\begin{figure}[H]
\centering
\includegraphics[width=0.9\textwidth]{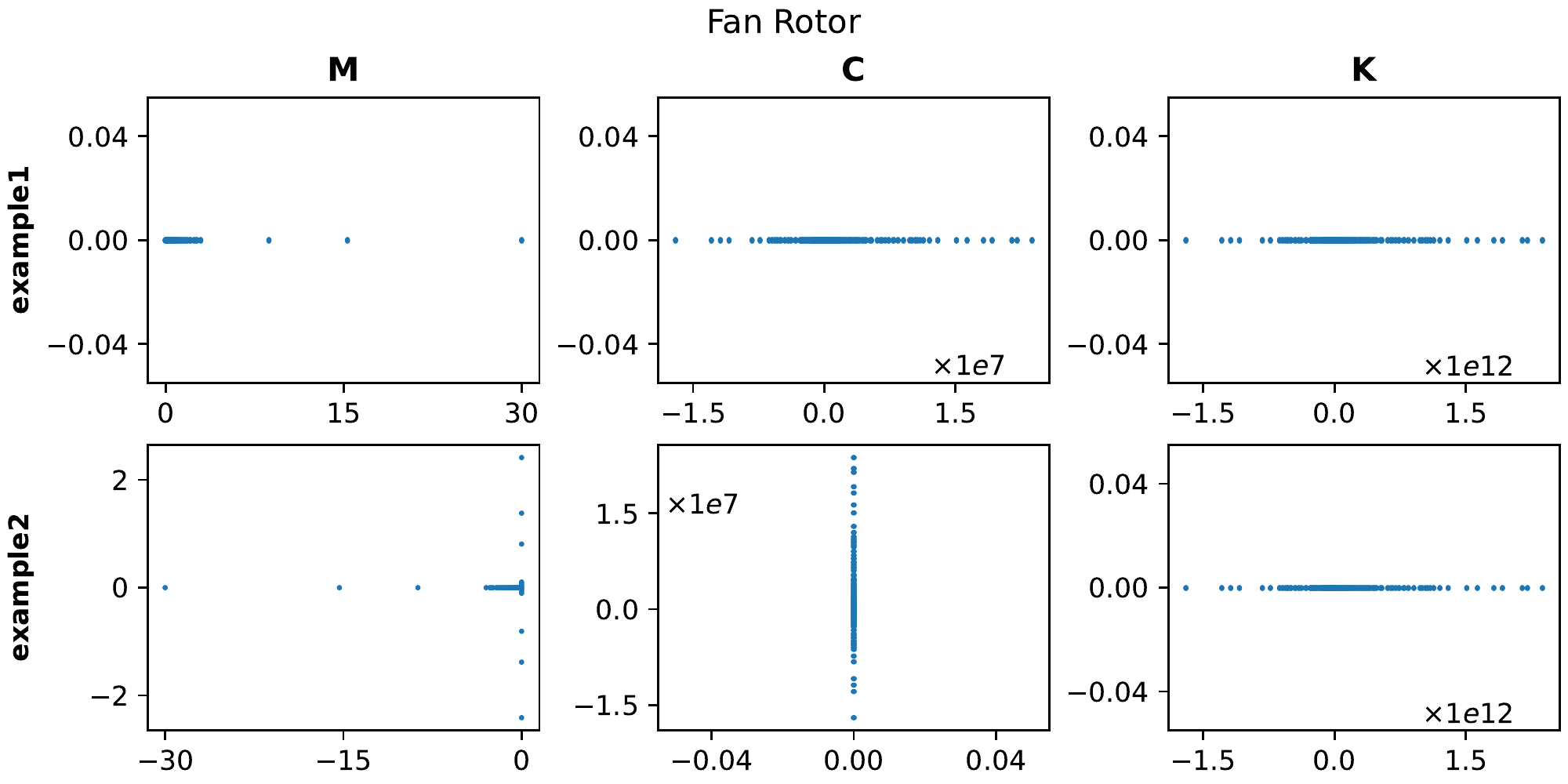}
\caption{The elements of the examples' M, C, K in the complex plane} 
\end{figure}

\end{appendix}

\bibliography{References}

\end{document}